\documentclass[aos,MSNbibl,seceqn,dvips]{arximspdf}
\usepackage{dcolumn}
\usepackage{graphicx}

\doi{10.1214/13-AOS1201} 
\volume{42}
\issue{2}
\pubyear{2014}
\firstpage{700}
\lastpage{727}

\makeatletter
\newcolumntype{d}[1]{D{.}{.}{#1}}
\newtheorem{teo}{Theorem}
\newtheorem{lem}{Lemma}
\newproclaim{rem}{Remark}
\newproclaim{exa}{Example}
\makeatother

\begin{document}
\begin{frontmatter}

\title{A semiparametric spatial dynamic model\thanksref{T1}}
\runtitle{SSDM}

\begin{aug}
\author[A]{\fnms{Yan} \snm{Sun}\ead[label=e1]{sunyan@mail.shufe.edu.cn}},
\author[B]{\fnms{Hongjia} \snm{Yan}\ead[label=e2]{yanli8626@qq.com}},
\author[B]{\fnms{Wenyang} \snm{Zhang}\corref{}\ead[label=e3]{wenyang.zhang@york.ac.uk}}
\and
\author[C]{\fnms{Zudi} \snm{Lu}\ead[label=e4]{Z.Lu@soton.ac.uk}}
\runauthor{Sun, Yan, Zhang and Lu}
\affiliation{Shanghai University of Finance and Economics, The
University of York,\\ The University of York and The University of Southampton}
\address[A]{Y. Sun\\
School of Economics\\
Shanghai University of Finance and Economics\\
No. 777 Guoding Road\\
Shanghai, 200433\\
P.R. China\\
\printead{e1}} 
\address[B]{H. Yan\\
W. Zhang\\
Department of Mathematics\\
The University of York\\
Heslington\\
York, YO10 5DD\\
United Kingdom\\
\printead{e2}\\
\phantom{E-mail:\ }\printead*{e3}}
\address[C]{Z. Lu\\
School of Mathematical Sciences\\
The University of Southampton\\
Highfield\\
Southampton, SO17 1BJ\\
United Kingdom\\
\printead{e4}}
\end{aug}
\thankstext{T1}{Supported by National Science Foundation of China
(Grant 11271242), program for New
Century Excellent Talents in China's University (NECT-10-0562),
Key Laboratory of Mathematical Economics (SUFE), Ministry of Education
of China
and the Singapore National Research Foundation under its
Cooperative Basic Research Grant and administered by the Singapore
Ministry of
Health's National Medical Research Council (Grant No. NMRC/CBRG/0014/2012).}

\received{\smonth{5} \syear{2013}}
\revised{\smonth{9} \syear{2013}}

%
\begin{abstract}
Stimulated by the Boston house price data, in this paper, we propose a
semiparametric spatial dynamic model, which extends the
ordinary spatial autoregressive models to accommodate the effects of some
covariates associated with the house price. A profile likelihood based
estimation procedure is proposed. The asymptotic normality of the proposed
estimators are derived. We also investigate how to identify the
parametric/nonparametric components in the proposed semiparametric model.
We show how many unknown parameters an unknown bivariate function
amounts to,
and propose an AIC/BIC of nonparametric version for model selection.
Simulation studies are
conducted to examine the performance of the proposed methods. The simulation
results show our methods work very well. We finally apply the proposed
methods to analyze the Boston house price data, which leads to some interesting
findings.
\end{abstract}

%

\end{frontmatter}

\section{Introduction}\label{sec1}

The Boston house price data is frequently used in literature to illustrate
some new statistical methods. If we use $y_i$ to denote the median
value of
owner-occupied homes at location $s_i$, a spatial autoregressive model for
the data would be
%
\begin{equation}
y_i = \sum_{j \neq i} w_{ij}
y_j + \varepsilon_i, \qquad i = 1, \ldots, n,
\label{moti0}
\end{equation}
where $w_{ij}$ is the impact of $y_j$ on $y_i$.
However, (\ref{moti0}) is inadequate because it models $y_i$ solely
based on the median value prices, $y_j,$ for $j \neq i$.
It is
better to
incorporate the effects of some important covariates, such as the crime
rate and
accessibility to radial highways, into the model. Let $X_i$, a $p$-dimensional
vector, be the vector of the covariates associated with $y_i$. A~reasonable
model to fit the data would be
%
\begin{equation}
y_i = \sum_{j \neq i} w_{ij}
y_j + X_i^{\mathrm{T}} \bolds{\beta}+
\varepsilon_i, \qquad i = 1, \ldots, n, \label{moti1}
\end{equation}
where $w_{ij}$ and $\bolds{\beta}$ are unknown. However, there are two
problems with
model
(\ref{moti1}): first, there are too many unknown parameters; second, the
model has not taken into account the location effects of the impacts of the
covariates---the impacts of some covariates may vary over location. To
control the number of unknown parameters and take the location effects into
account, we propose the following model to fit the data:
%
\begin{equation}
y_i = \alpha\sum_{j \neq i} w_{ij}
y_j + X_i^{\mathrm{T}} \bolds{\beta}(s_i)
+ \varepsilon_i, \qquad i = 1, \ldots, n, \label{model}
\end{equation}
where $w_{ij}$ is a specified certain physical or economic distance,
$s_i$ is
the location of the $i$th observation, which is a two-dimensional vector,
$\bolds{\beta}(\cdot) = (\beta_1(\cdot), \ldots, \beta
_p(\cdot
))^{\mathrm{T}}$,
$\varepsilon_i$, $i = 1, \ldots, n$, are i.i.d., and follow
$N(0, \sigma^2)$, $\{X_i, i = 1, \ldots, n \}$ is
independent of $\{\varepsilon_i, i = 1, \ldots, n \}$. $\alpha$,
$\sigma^2$ and $\bolds{\beta}(\cdot)$ are unknown to
be estimated. Model~(\ref{model}) is the model this paper is going to
address. From now on, $y_i$ is of course not necessarily the house
price, it
is a generic response variable. We will also see that the normality assumption
imposed on $\varepsilon_i$ is just for the description of the
construction of the
proposed estimation procedure. It is not necessary for the asymptotic
properties of the proposed estimators.

In model~(\ref{model}), the spatial neighboring effect of $y_j$, $j
\neq i$,
on $y_i$ is formulated through $\alpha w_{ij}$, where $w_{ij}$ is a specified
certain physical or economic distance, and $\alpha$ is an unknown
baseline of
the spatial neighboring effect. Such method to define spatial neighboring
effect is common; see Ord \cite{r13}, Anselin \cite{r1}, Su and Jin
\cite{r14}.

If there is no any condition imposed on the spatial neighboring
effects, and
the spatial neighboring effects are formulated as unknown $w_{ij}$,
$i=1, \ldots, n$, \mbox{$j = 1, \ldots, i-1, i+1, \ldots, n$}, we would
have $(n - 1)n$ unknown $w_{ij}$'s to estimate. In which case, it would be
impossible to have consistent estimators of $w_{ij}$'s. However, if we impose
some kind of sparsity on $w_{ij}$'s, by penalized maximum likelihood estimation,
it is possible to construct consistent estimators of $w_{ij}$'s.
However, that
has gone beyond the scope of this paper although it is a promising research
project.

Model~(\ref{model}) is a useful extension of spatial autoregressive models
(Gao et al.~\cite{r6}; Kelejian and Prucha \cite{r9}; Ord \cite{r13};
Su and Jin \cite{r14}) and
varying coefficient models (Cheng et~al. \cite{r2}; Fan and Zhang \cite
{r4,r5};
Li and Zhang \cite{r11}; Sun et~al. \cite{r15}; Zhang et~al. \cite{r19,r20}; Wang and Xia \cite{r17};
and Tao and Xia \cite{r16}). One characteristic of model
(\ref{model}) is
\[
E(\varepsilon_i|y_1, \ldots, y_{i-1},
y_{i+1}, \ldots, y_n) \neq0
\]
although $E(\varepsilon_i) = 0$, the standard least squares estimation
will not
work for~(\ref{model}). In this paper, based on the local linear
modeling and profile likelihood idea, we will propose a local likelihood
based
estimation procedure for the unknown parameters and functions in~(\ref{model})
and derive the asymptotic properties of the obtained estimators.

In reality, some components of $\bolds{\beta}(\cdot)$ in model (\ref
{model}) may be
constant, and we do not know which components are functional, which are
constant. Methodologically speaking, if mistakenly treating a constant
component as functional, we would pay a price on the variance side of the
obtained estimator; on the other hand, if mistakenly treating a functional
component as constant, we would pay a price on the bias side of the obtained
estimator. The identification of constant/functional components in
$\bolds{\beta}(\cdot)$ is imperative. From practical point of view, the
identification of constant components is also of importance. For the
data set
we study in this paper, $\bolds{\beta}(\cdot)$ can be interpreted as the
vector of the
impacts of the covariates concerned on the house price. The identification
will
reveal which covariates have location varying impacts on the house
price, and which
do not. This is apparently something of great interest. In this paper,
we will show how many unknown parameters an unknown bivariate function amounts
to, and propose an AIC/BIC of nonparametric version to identify the constant
components of $\bolds{\beta}(\cdot)$ in model~(\ref{model}).

The paper is organized as follows. We begin in Section~\ref{sec2} with a
description of
the estimation procedure for the proposed model~(\ref{model}). In
Section~\ref{sec3}, we
show how many unknown parameters an unknown bivariate function amounts
to, and
propose an AIC/BIC of nonparametric version for model selection. Asymptotic
properties of the proposed estimators are presented in Section~\ref{sec4}. The
performance of the proposed methods, including both estimation and model
selection methods, is assessed by a simulation study in Section~\ref{sec5}. In
Section~\ref{sec6}, we explore how the covariates, which are commonly found to
be associated
with house price, affect the median value of owner-occupied homes in Boston,
and how the impacts of these covariates change over location based on the
proposed model and estimation procedure.

Throughout this paper, $\mathbf{0}_k$ is a $k$-dimensional vector with each
component being $0$, $I_k$ is an identity matrix of size $k$, $U[0,
1]^2$ is
a two-dimensional uniform distribution on $[0, 1] \times[0, 1]$.

\section{Estimation procedure}\label{sec2}

Let $w_{ii} = 0$, $W = (w_{ij})$, $Y = (y_1, \ldots, y_n)^{\mathrm{T}}$,
$A = I_n - \alpha W$ and
$
\mathbf{m}
=
(X_1^{\mathrm{T}} \bolds{\beta}(s_1), \ldots, X_n^{\mathrm{T}} \bolds
{\beta}(s_n)
)^{\mathrm{T}}
$. By simple calculations, we have that the conditional density
function of $Y$\vadjust{\goodbreak}
given $\mathbf{m}$ is
$N (
A^{-1} \mathbf{m},\break  (A^{\mathrm{T}} A)^{-1} \sigma^2
)$,
which leads to the following log likelihood function
%
\begin{equation}
-\frac{n}{2} \log(2 \pi) - n \log(\sigma) + \log\bigl(|A|\bigr) -\frac{1}{2\sigma^2}(A
Y - \mathbf{m})^{\mathrm{T}} (A Y - \mathbf{m}). \label{lik}
\end{equation}

Our estimation is profile likelihood based. We first construct the estimator
$\tilde{\bolds{\beta}}(\cdot; \alpha)$ of $\bolds{\beta}(\cdot
)$ pretending
$\alpha$ is
known, then let $(\hat{\alpha}, \hat{\sigma}^2)$ maximize~(\ref{lik})
with $\bolds{\beta}(\cdot)$ being replaced by $\tilde{\bolds{\beta
}}(\cdot; \alpha)$.
$\hat{\alpha}$ and $\hat{\sigma}^2$ are our estimators of $\alpha$ and
$\sigma^2$, respectively. After the estimator of $\alpha$ is
obtained, the
estimator of $\bolds{\beta}(\cdot)$ is taken to be $\tilde{\bolds{\beta
}}(\cdot;
\alpha)$
with $\alpha$ and the bandwidth used being replaced by $\hat{\alpha
}$ and a
slightly larger bandwidth, respectively. The details are as follows.

For\vspace*{1.5pt} any $s=(u, v)^{\mathrm{T}}$, we denote
$(\partial\bolds{\beta}(s) / \partial u, \partial\bolds{\beta
}(s) / \partial
v)$ by
$\dot{\bolds{\beta}}(s)$, where
$\partial\bolds{\beta}(s) / \partial u
=
(\partial\beta_1(s) / \partial u, \ldots, \partial\beta_p(s) / \partial
u)^{\mathrm{T}}$. We define $\|s\| =
(s^{\mathrm{T}}
s)^{1/2}$.

For any given $s$, by the Taylor's expansion, we have
\[
\bolds{\beta}(s_i) \approx\bolds{\beta}(s) + \dot{\bolds{\beta
}}(s) (s_i - s),
\]
when $s_i$ is in a small neighborhood of $s$, which leads to the following
objective function for estimating $\bolds{\beta}(s)$:
%
\begin{equation}
\sum_{i=1}^n \bigl( y_i^* -
X_i^{\mathrm{T}} \mathbf{a}- X_i^{\mathrm{T}}
\mathbf{B}(s_i - s) \bigr)^2 K_h\bigl(\|s_i - s\|\bigr), \label{loc}
\end{equation}
where $y_i^*$ is the $i$th component of $AY$, $K_h(\cdot) = K(\cdot/h)/h^2$,
$K(\cdot)$ is a kernel function, and $h$ is a bandwidth.
Let $(\hat{\mathbf{a}}, \widehat{\mathbf{B}})$
minimise~(\ref{loc}), the ``estimator'' $\tilde{\bolds{\beta}}(s; \alpha
)$ of
$\bolds{\beta}(s)$ is taken to be $\hat{\mathbf{a}}$. By simple
calculations, we have
%
\begin{equation}
\tilde{\bolds{\beta}}(s; \alpha) = \hat{\mathbf{a}} = (I_p,
\mathbf{0}_{p \times2p}) \bigl(\mathcal{X}^{\mathrm{T}} \mathcal {W}
\mathcal{X} \bigr)^{-1} \mathcal{X}^{\mathrm{T}} \mathcal{W}A Y,
\label{est}
\end{equation}
where $\mathbf{0}_{p \times q}$ is a matrix of size $p \times q$ with each
entry being
$0$, and
\begin{eqnarray*}
\mathcal{X} &=& \pmatrix{ X_1 & \cdots& X_n
\vspace*{3pt}\cr
X_1 \otimes(s_1 - s) & \cdots& X_n
\otimes(s_n - s) }^{\mathrm{T}},
\\
\mathcal{W} &=& \operatorname{diag} \bigl( K_h\bigl(\|s_1 - s \|\bigr), \ldots, K_h\bigl(\|s_n - s\|\bigr) \bigr).
\end{eqnarray*}
Replacing $\bolds{\beta}(s_i)$ in~(\ref{lik}) by $\tilde{\bolds{\beta
}}(s_i; \alpha)$
and ignoring the constant term, we have the objective function for estimating
$\alpha$ and $\sigma^2$
%
\begin{equation}
- n \log(\sigma) + \log\bigl(|A|\bigr) -\frac{1}{2\sigma^2}(A Y - \tilde{
\mathbf{m}})^{\mathrm{T}} (A Y - \tilde{\mathbf{m}}), \label{profil}
\end{equation}
where $\tilde{\mathbf{m}}$ is $\mathbf{m}$ with $\bolds{\beta
}(s_i)$ being
replaced by
$\tilde{\bolds{\beta}}(s_i; \alpha)$. Let $\alpha_i$, $i=1, \ldots, n$, be
the eigenvalues of $W$,
\[
\tilde{\sigma}^2 = \frac{1}{n} (A Y - \tilde{
\mathbf{m}})^{\mathrm{T}} (A Y - \tilde{\mathbf{m}})
\]
and $(\hat{\alpha}, \hat{\sigma}^2)$ maximize~(\ref{profil}).
Noticing that
$|A| = \prod_{i=1}^n(1 - \alpha\alpha_i)$, by simple
calculations, we
have $\hat{\alpha}$ is the maximizer of
%
\begin{equation}
- n \log(\tilde{\sigma}) + \sum_{i=1}^n
\log \bigl(|1- \alpha\alpha_i| \bigr) \label{lam}
\end{equation}
and $\hat{\sigma}^2$ is $\tilde{\sigma}^2$ with $\alpha$ being
replaced by
$\hat{\alpha}$.

Note that the maximization of (\ref{lam}) is not difficult because it is
a one-dimensional optimization problem, which can be solved using a grid point
method.

The estimator
$\hat{\bolds{\beta}}(\cdot)$
$ (
=
(\hat{\beta}_1(\cdot), \ldots, \hat{\beta}_p(\cdot
))^{\mathrm{T}}
)$
is $\tilde{\bolds{\beta}}(\cdot; \alpha)$
with $\alpha$ being replaced by~$\hat{\alpha}$ and the bandwidth $h$
by a
slightly larger bandwidth $h_1$.
The reason for replacing the bandwidth $h$ by a slightly larger number $h_1$
is that the former bandwidth is appropriate for the estimation of
constant parameters,
$\alpha$ and $\alpha^2$, and the latter is more appropriate for
the estimation of functional parameters. Also, the estimators
of constant parameters need a smaller bandwidth $h$ in order
to achieve the optimal rate of convergence.

In reality, some components of $\bolds{\beta}(\cdot)$ may be
constant. If a
component of $\bolds{\beta}(\cdot)$ is a constant, say $\beta
_1(\cdot) =
\beta_1$, we
use the average of $\hat{\beta}_1(s_i)$, $i=1, \ldots, n$, to estimate
the constant $\beta_1$, that is,
\[
\hat{\beta}_1 = \frac{1}{n} \sum
_{i=1}^n \hat{\beta}_1(s_i).
\]
How to identify the constant components of $\bolds{\beta}(\cdot)$
will be addressed
in the next section.

\section{Identification of constant components}\label{sec3}
\subsection{Criterion for identification}\label{sec3.1}

As we mentioned before, some components of $\bolds{\beta}(\cdot)$ in model
(\ref{model}) may be constant in reality, and to identify such constant
components is of importance. In this paper, we appeal the $\mathrm{AIC}$ or
$\mathrm{BIC}$ to
identify the constant components. The $\mathrm{AIC}$ for (\ref
{model}), in
which some
components of $\bolds{\beta}(\cdot)$ may be constant, is defined as follows:
%
\begin{equation}
\mathrm{AIC} = n \log(\hat{\sigma}) - \log\bigl(|\widehat{A}|\bigr) + \frac
{1}{2\hat{\sigma}^2}(
\widehat{A} Y - \hat{\mathbf{m}})^{\mathrm{T}} (\widehat{A} Y - \hat{
\mathbf{m}}) + \mathcal{K}, \label{aic}
\end{equation}
where $\widehat{A}$ and $\hat{\mathbf{m}}$ are $A$ and $\mathbf{m}$
with the unknown
parameters and functions being replaced by their estimators, $\mathcal{K}$
is the
number of unknown parameters in model~(\ref{model}). The $\mathrm
{BIC}$ can be
defined in a similar way.

Because there are unknown functions in model~(\ref{model}), the first hurdle
in the calculation of $\mathrm{AIC}$ of model~(\ref{model}) is to
find how
many unknown
constants an unknown bivariate function amounts to. In the following, based
on the residual sum of squares of standard bivariate nonparametric regression
model, we propose an ad hoc way to solve this problem.

Suppose we have the following standard bivariate nonparametric regression
model:
%
\begin{equation}
\eta_i = g(s_i) + e_i, \qquad i = 1,
\ldots, n, \label{biv}
\end{equation}
where $E(e_i)= 0$ and $\operatorname{var}(e_i) = \sigma_e^2$. The
residual sum of squares
of~(\ref{biv}) is
\[
\mathrm{RSS}= \sum_{i=1}^n \bigl\{
\eta_i - \hat{g}(s_i) \bigr\}^2,
\]
where $\hat{g}(\cdot)$ is the local linear estimator of $g(\cdot)$.
On the
other hand,
\[
\hspace*{-5pt}E \bigl(\mathrm{RSS}/\sigma_e^2
\bigr) = n-\mbox{the number of unknown parameters in the regression function.}
\]
So, the number $\mathcal{T}$ of unknown constants the unknown function
$g(\cdot)$
amounts to can be reasonably viewed as
\[
\mathcal{T} = n - E \bigl(\mathrm{RSS}/\sigma_e^2
\bigr) = n - \sigma_e^{-2} E \Biggl[ \sum
_{i=1}^n \bigl\{\eta_i -
\hat{g}(s_i) \bigr\}^2 \Biggr].
\]
To make $\mathcal{T}$ more convenient to use, we derive the asymptotic
form of
$\mathcal{T}$. Let
\[
\mathbf{S}_i = \pmatrix{1 & s_1^{\mathrm{T}} -
s_i^{\mathrm{T}}
\vspace*{3pt}\cr
\vdots& \vdots
\vspace*{3pt}\cr
1 & s_n^{\mathrm{T}} - s_i^{\mathrm{T}}}, \qquad
\bolds{\eta} = \pmatrix{ \eta_1
\vspace*{3pt}\cr
\vdots
\vspace*{3pt}\cr
\eta_n},
\qquad\mathbf{e} = \pmatrix{e_1
\vspace*{3pt}\cr
\vdots
\vspace*{3pt}\cr
e_n}
\]
and
\[
\mathcal{W}_i = \operatorname{diag} \bigl( K_h(u_1
- u_i)K_h(v_1 - v_i), \ldots,
K_h(u_n - u_i)K_h(v_n
- v_i) \bigr),
\]
we have
\[
\hat{g}(s_i) = (1, 0, 0) \bigl( \mathbf{S}_i^{\mathrm{T}}
\mathcal{W}_i \mathbf{S}_i \bigr)^{-1}
\mathbf{S}_i^{\mathrm{T}} \mathcal{W}_i \bolds{\eta}.
\]
By the standard argument in Fan and Gijbels \cite{r3} and Lemma~1 in
Fan and Zhang~\cite{r4}, we have
\[
\mathcal{T}= \bigl(2 K^2(0) - \nu_*^2 \bigr)
h^{-2} + o \bigl(h^{-2} \bigr),
\]
when $h = o(n^{-1/6})$ and $n h^2 \rightarrow\infty$, where
$\nu_* = \int K^2(t) \,dt$.

We conclude that an unknown bivariate function amounts to
$ (2 K^2(0) - \nu_*^2 ) h^{-2}$ unknown constants. Based on this
conclusion, if the number of constant components in $\bolds{\beta
}(\cdot)$
is $q$,
the $\mathcal{K}$ in~(\ref{aic}) will be
$q + (p - q) (2 K^2(0) - \nu_*^2 ) h^{-2}$.

To identify the constant components in $\bolds{\beta}(\cdot)$ in
(\ref
{model}) is
basically a model selection problem. Theoretically speaking, we go for the
model with the smallest $\mathrm{AIC}$ (or $\mathrm{BIC}$). However, in
practice, it is
almost computationally impossible to compute the AICs for all possible models.
We have to use some algorithm to reduce the computational burden. In the
following, we are going to introduce two algorithms for the model selection.

\subsection{Computational algorithms}\label{sec3.2}

In this section, we use $\mathrm{AIC}$ as an example to demonstrate
the introduced
algorithms. The model in which $\bolds{\beta}(\cdot)$ has its
$i_1$th, $i_2$th$,\ldots,i_k$th
components being constant is denoted by
$\{i_1, \ldots, i_k\}$. When $k=0$, we define the model as the
model in
which all components of $\bolds{\beta}(\cdot)$ are functional, and
denote it by
$\{\mbox{ }\}$.

\subsubsection*{Backward elimination}
The first algorithm we introduce is the backward elimination. Details
are as
follows.
\begin{longlist}[(2)]
\item[(1)] We start with the full model, $\{1, \ldots, p\}$, and compute
its $\mathrm{AIC}$ by~(\ref{aic}). Denote the full model by $\mathcal
{M}_p$, its
$\mathrm{AIC}$ by
$\mathrm{AIC}_p$.

\item[(2)] For any integer $k$, suppose the current model is
$\mathcal{M}_k = \{i_1, \ldots, i_k\}$ with $\mathrm{AIC}$ given by
$\mathrm{AIC}_k$. Take
$\mathcal{M}_{k-1}$ to be the model with the largest maximum of log likelihood
function among the models
$\{i_1, \ldots, i_{j-1}, i_{j+1}, \ldots, i_k \}$,
$j=1, \ldots, k$. If $\mathrm{AIC}_k < \mathrm{AIC}_{k-1}$, the
chosen model
is $\mathcal{M}_k$, and the model selection is ended; otherwise,
continue to compute
$\mathcal{M}_l$ and $\mathrm{AIC}_l$ until either $\mathrm{AIC}_l <
\mathrm{AIC}_{l-1}$ or $l=0$.
\end{longlist}

\subsubsection*{Curvature-to-average ratio (CTAR) based method}

A more aggressive way to reduce the computational burden involved in
the model
selection procedure is based on the ratio of the curvature of the estimated
function to its average. Explicitly, we first treat all $\beta_j(\cdot)$,
$j=1, \ldots, p$, as functional. For each $j$, $j=1, \ldots,
p$, we compute the curvature-to-average ratio (CTAR) $R_j$ of the estimated
function~$\hat{\beta}_j(\cdot)$:
\[
R_j = \frac{1}{\bar{\beta}_j^2}\sum_{i=1}^n
\bigl\{ \hat{\beta}_j(s_i) - \bar{\beta}_j
\bigr\}^2, \qquad\bar{\beta}_j = \frac{1}{n} \sum
_{i=1}^n \hat{\beta}_j(s_i),
\qquad j = 1, \ldots, p.
\]
We sort $R_j$, $j=1, \ldots, p$, in an increasing order, say
$
R_{i_1} \leq\cdots\leq R_{i_p}$,
then compute the $\mathrm{AIC}$s for the models $\{i_1, \ldots, i_k\}$
from $k=0$
to the turning point $k_0$ where the $\mathrm{AIC}$ starts to
increase. The chosen
model is $\{i_1, \ldots, i_{k_0}\}$.

The algorithm based on the CTAR is much faster than the backward elimination
based algorithm, however, we find it less accurate although it still works
reasonably well in our simulation studies. This is because the CTARs of all
coefficients are obtained in one go based on the model in which all
coefficients are treated as functional, and not updated. This will
speed up
the selection procedure; on the other hand, the effect of randomness
would be
stronger than that in backward elimination, which leads to a slightly larger
possibility of picking up a wrong model.

\section{Asymptotic properties}\label{sec4}

In this section, we are going to present the asymptotic properties of the
proposed estimators. We will, in this section, only present the asymptotic
results, and leave the theoretical proofs in the \hyperref[app]{Appendix}.

Although we assume $\varepsilon_i$ in~(\ref{model}) follows normal distribution
in our model assumption, we do not need this assumption when deriving the
asymptotic properties of the proposed estimators. So, in this section,
we do
not assume $\varepsilon_i$ follows normal distribution unless otherwise stated.

In this section, for $w_{ij}$ in~(\ref{model}), we assume that there exists
a sequence $\rho_n>0$ such that $w_{ij} = O(1/\rho_n)$ uniformly with respect
to $i, j$ and the matrices $W$ and $A^{-1}$ are uniformly bounded in
both row
and column sums.

We now introduce some notations needed in the presentation of the asymptotic
properties of the proposed estimators: let $\mu_j=E\varepsilon_1^j,
j=1,\ldots,4$,
\begin{eqnarray*}
\kappa_0&=&\int_{R^2} K\bigl(\|s\|\bigr) \,ds,
\\
\kappa_2&=&\int_{R^2} \bigl[(1,0)s \bigr]^2 K\bigl(\|s\|\bigr) \,ds =\int_{R^2} \bigl[(0,1)s
\bigr]^2 K\bigl(\|s\|\bigr) \,ds,
\\
\nu_0&=&\int_{R^2} K^2\bigl(\|s\|\bigr) \,ds,
\\
\nu_2&=&\int_{R^2} \bigl[(1,0)s \bigr]^2
K^2\bigl(\|s\|\bigr) \,ds= \int_{R^2} \bigl[(0,1)s
\bigr]^2 K^2\bigl(\|s\|\bigr) \,ds,
\\
G &=& (g_{ij}) = W A^{-1},\qquad \Psi= E \bigl(X_1X_1^{\mathrm{T}}\bigr),\qquad \Gamma=EX_1,
\\
Z_1(s)&=&\lim_{n\rightarrow\infty} \frac{1}{n } \sum
_{i=1}^n g_{ii}\bolds{
\beta}(s_i) K_h\bigl(\|s_i-s\|\bigr),
\\
Z_2(s)&=&\lim_{n\rightarrow\infty} \frac{1}{n } \sum
_{i=1}^n\sum_{j\neq i}^n
g_{ij} \bolds{\beta}(s_j)K_h\bigl(\|s_i-s\|\bigr),
\\
Z(s) &=& Z_1(s)+\Psi^{-1}\Gamma\Gamma^{\mathrm{T}}Z_2(s),
\\
Z&=&\kappa_0^{-1} \bigl(f^{-1}(s_1)X_1^{\mathrm{T}}Z(s_1),
\ldots,f^{-1}(s_n)X_n^{\mathrm{T}}Z(s_n)
\bigr)^{\mathrm{T}},
\\
\pi_1&=&\lim_{n\rightarrow\infty} \frac{\operatorname{tr}((G+G^{\mathrm{T}})G)}{
n}, \qquad
\pi_2=\lim_{n\rightarrow\infty} \frac{\operatorname
{tr}(G)}{ n},
\\
\pi_3&=&\lim_{n\rightarrow\infty} \frac{1}{n} \sum
_{i=1}^n g_{ii}^2,
\\
\lambda_{1}&=&\lim_{n\rightarrow\infty} \frac{1}{n}E
\bigl[(G\mathbf {m}-Z)^{\mathrm{T}}(G \mathbf{m} - Z) \bigr],
\\
\lambda_{2}&=&\lim_{n\rightarrow\infty} \frac{1}{n}E
\bigl[(G\mathbf{m}-Z)^{\mathrm{T}}G_c \bigr], \qquad
\lambda_{3}=\lim_{n\rightarrow\infty} \frac{1}{n}E \bigl[(G
\mathbf{m}-Z)^{\mathrm{T}}\mathbf{1}_n \bigr],
\end{eqnarray*}
where $G_c=(g_{11},\ldots, g_{nn})^{\mathrm{T}}$ and $\mathbf{1}_n$
is an $n$-dimensional vector
with each component being $1$. Further, let
\begin{eqnarray*}
\Omega&=& \pmatrix{\displaystyle\frac{1 } {\sigma^2}\lambda_1+\pi_1
&\displaystyle \frac{ 1}{\sigma^2}\pi_2
\vspace*{6pt}\cr
\displaystyle\frac{1 }{\sigma^2}\pi_2
&\displaystyle \frac{1}{2\sigma^4}},
\\[3pt]
\Sigma&=&\pmatrix{\displaystyle\frac{\mu_4-3\sigma^4}{\sigma^4}\pi_3 +\frac{2\mu_3}{\sigma^4}\lambda_2
&\displaystyle \frac{\mu_3}{2\sigma^6} \lambda_3+ \frac{\mu_4-3\sigma^4}{2\sigma^6}\pi_2
\vspace*{6pt}\cr
\displaystyle\frac{\mu_3}{2\sigma^6}\lambda_3+ \frac{\mu_4-3\sigma^4}{2\sigma^6}\pi_2
&\displaystyle \frac{\mu_4-3\sigma^4}{4\sigma^8}},
\\
s&=&(u,v)^{\mathrm{T}}, \qquad\bolds{\beta}_{uu}(s)= \biggl(
\frac{\partial^2 \beta_1(s)}{\partial
u^2},\ldots, \frac{\partial^2 \beta_p(s)}{\partial u^2} \biggr)^{\mathrm{T}},
\\
\bolds{\beta}_{vv}(s)&=& \biggl(\frac{\partial^2 \beta_1(s)}{\partial
v^2},\ldots,
\frac{\partial^2 \beta_p(s)}{\partial v^2} \biggr)^{\mathrm{T}}
\end{eqnarray*}
and
\[
S = \pmatrix{ \bigl(X_1^{\mathrm{T}}, \mathbf{0}_{1\times2p}
\bigr) \bigl(\mathcal{X}_{(1)}^{\mathrm{T}}\mathcal{W}_{(1)}
\mathcal{X}_{(1)} \bigr)^{-1} \mathcal{X}_{(1)}^{\mathrm{T}}
\mathcal{W}_{(1)}
\vspace*{3pt}\cr
\vdots
\vspace*{3pt}\cr
\bigl(X_n^{\mathrm{T}}, \mathbf{0}_{1\times2p} \bigr) \bigl(
\mathcal{X}_{(n)}^{\mathrm{T}}\mathcal{W}_{(n)}\mathcal
{X}_{(n)} \bigr)^{-1} \mathcal{X}_{(n)}^{\mathrm{T}}
\mathcal{W}_{(n)}},
\]
where $\mathcal{X}_{(i)}$ and $\mathcal{W}_{(i)}$ are $\mathcal{X}$
and $\mathcal{W}$, respectively,
with $s$ being replaced by \mbox{$s_i$}, \mbox{$i=1, \ldots, n$}.

By some simple calculations, we can see the matrix $\Omega$ defined
above is
the limit of the Fisher information matrix of $\alpha$ and $\sigma
^2$. As the
singularity of matrix $\Omega$ may have serious implication on the
convergence rate of the proposed estimators, we present the asymptotic
properties for the case where $\Omega$ is nonsingular and the case where
$\Omega$ is singular separately. We present the nonsingular case in Theorems
\ref{teo1}--\ref{teo3}, and singular case in Theorems~\ref{teo4}--\ref{teo7}.

\begin{teo}\label{teo1}
Under\vspace*{1pt} the conditions (1)--(7) or conditions (1)--(6),
($\tilde{7}$)\break  and (8) in \hyperref[app]{Appendix}, $\Omega$ is nonsingular, and when
$n^{1/2}h^2/\log^2 n\rightarrow\infty$ and $nh^8\rightarrow0$,
$\hat{\alpha}$~and~$\hat{\sigma}^2$ are consistent estimators of
$\alpha$~and~$\sigma^2$, respectively.
\end{teo}

Theorem~\ref{teo1} shows the conditions under which $\Omega$ is nonsingular and the
consistency of $\hat{\alpha}$ and $\hat{\sigma}^2$ under such conditions.
Based on Theorem~\ref{teo1}, we can derive the asymptotic normality of $\hat
{\alpha}$ and
$\hat{\sigma}^2$.

\begin{teo}\label{teo2}
Under the assumptions of Theorem~\ref{teo1}, if the second partial
derivative of $\bolds{\beta}(s)$ is Lipschitz continuous and
$nh^6\rightarrow0$,
\[
\sqrt{n} \bigl(\hat{\alpha} -\alpha, \hat{\sigma}^2-
\sigma^2 \bigr)^{\mathrm{T}} \stackrel{D} {\longrightarrow} N \bigl(
\mathbf{0}, \Omega^{-1} + \Omega^{-1}\Sigma
\Omega^{-1} \bigr).
\]
Further, if $\varepsilon_i$ is normally distributed,
\[
\sqrt{n} \bigl(\hat{\alpha} -\alpha, \hat{\sigma}^2-
\sigma^2 \bigr)^{\mathrm{T}} \stackrel{D} {\longrightarrow} N \bigl(
\mathbf{0}, \Omega^{-1} \bigr).
\]
\end{teo}

Theorem~\ref{teo2} implies that the convergence rate of $\hat{\alpha}$ is of order
$n^{-1/2}$ when $\Omega$ is nonsingular, which is the optimal rate for
parametric estimation. We will see, in Theorem~\ref{teo5}, this rate can not be
achieved by $\hat{\alpha}$ when $\Omega$ is singular.

\begin{teo}\label{teo3}
Under the assumptions of Theorem~\ref{teo1}, if
$nh_1^6=O(1)$ and
$h/h_1\rightarrow0$,
\begin{eqnarray*}
&& \sqrt{nh_1^2 f(s)} \bigl(\hat{\bolds{
\beta}}(s)-\bolds{\beta}(s)-2^{-1}\kappa_0^{-1}
\kappa_2 h_1^2 \bigl\{\bolds{
\beta}_{uu}(s)+ \bolds{\beta}_{vv}(s) \bigr\} \bigr)
\\
&&\qquad \stackrel{D} {\longrightarrow} N \bigl(\mathbf{0}, \kappa_0^{-2}
\nu_0 \sigma^2 \Psi^{-1} \bigr)
\end{eqnarray*}
for any given $s$.
\end{teo}

Theorem~\ref{teo3} shows $\hat{\bolds{\beta}}(\cdot)$ is asymptotic normal and
achieves the
convergence rate of order $n^{-1/6}$, which is the optimal rate for bivariate
nonparametric estimation.

We now turn to the case where $\Omega$ is singular.

\begin{teo}\label{teo4}
Under the conditions (1)--(6) and (9) in the
\hyperref[app]{Appendix}, $\Omega$ is singular, and if
$nh^8\rightarrow0$, $n^{1/2}h^2/\log^2 n\rightarrow\infty$,
$\rho_n\rightarrow\infty$,
$\rho_n h^4\rightarrow0$ and $nh^2/\rho_n \rightarrow
\infty$, $\hat{\alpha}$ is a consistent estimator of $\alpha$.
\end{teo}

\begin{teo}\label{teo5}
Under the assumptions of Theorem~\ref{teo4}, if the
second partial derivative of $\bolds{\beta}(s)$ is Lipschitz
continuous and
$nh^6\rightarrow0$,
\[
\sqrt{n/\rho_n} (\hat{\alpha} - \alpha) \stackrel{D} {
\longrightarrow}N \bigl(0, \sigma^2\lambda_4^{-1}
\bigr),
\]
where
\[
\lambda_4=\lim_{n\rightarrow\infty} \frac{\rho_n}{n} E \bigl[(G
\mathbf{m} -SG \mathbf{m})^{\mathrm
{T}}(G\mathbf{m}-SG\mathbf{m}) \bigr].
\]
\end{teo}

Theorem~\ref{teo5} shows the convergence rate of $\hat{\alpha}$ is of order
$(n/\rho_n)^{-1/2}$ which is slower than $n^{-1/2}$ when
$\rho_n \rightarrow\infty$. However, we will see, from Theorem~\ref{teo7}, this
has no effect on the asymptotic properties of $\hat{\bolds{\beta
}}(\cdot)$.

\begin{teo}\label{teo6}
Under the assumptions of Theorem~\ref{teo5},
\[
\sqrt{n} \bigl(\hat{\sigma}^2 -\sigma^2 \bigr)
\stackrel{D} {\longrightarrow}N \bigl(0, \mu_4-\sigma^4
\bigr).
\]
\end{teo}

Theorem~\ref{teo6} shows that although the asymptotic variance of $\hat{\sigma
}^2$ is
different to that when $\Omega$ is nonsingular, $\hat{\sigma}^2$
still enjoys
convergence rate of $n^{-1/2}$.

\begin{teo}\label{teo7}
Under the assumptions of Theorem~\ref{teo4}, if
$nh_1^6=O(1)$ and $h/h_1\rightarrow0$,
\begin{eqnarray*}
&& \sqrt{nh_1^2 f(s)} \bigl(\hat{\bolds{
\beta}}(s)-\bolds{\beta}(s)- 2^{-1}\kappa_0^{-1}
\kappa_2 h_1^2 \bigl\{ \bolds{
\beta}_{uu}(s)+\bolds{\beta}_{vv}(s) \bigr\} \bigr)
\\
&&\qquad \stackrel{D} {\longrightarrow}N \bigl(\mathbf{0}, \kappa_0^{-2}
\nu_0 \sigma^2 \Psi^{-1} \bigr)
\end{eqnarray*}
for any given $s$.
\end{teo}

From Theorems~\ref{teo3} and~\ref{teo7}, we can see the singularity of $\Omega$
has no
effect on the asymptotic distribution of $\hat{\bolds{\beta}}(\cdot)$.

\section{Simulation studies}\label{sec5}

In this section, we will use simulated examples to examine the
performances of
the proposed estimation and model selection procedure. In all simulated
examples and the real data analysis later on, we set $w_{ij}$ to be
%
\begin{equation}
w_{ij} = \exp\bigl(-\|s_i - s_j\|\bigr) \bigm/ \sum
_{k \neq i} \exp\bigl(-\|s_i - s_k\|\bigr).
\label{weight}
\end{equation}

We first examine the performance of the proposed estimation procedure, then
the model selection procedure.

\subsection{Performance of the estimation procedure}\label{sec5.1}

\begin{exa}\label{exa1} In model~(\ref{model}), we set $p = 3$, $\sigma^2 = 1$,
\begin{eqnarray*}
\alpha&=& 0.5, \qquad\beta_1(s) = \sin \bigl(\|s\|^2 \pi
\bigr),\qquad
\beta_2(s) = \cos \bigl(\|s\|^2 \pi \bigr),
\qquad\beta_3(s) = e^{(\|s\|^2)}
\end{eqnarray*}
and independently generate $X_i$ from $N(\mathbf{0}_3, I_3)$, $s_i$ from
$U[0, 1]^2$, $\varepsilon_i$ from $N(0, \sigma^2)$, $i=1, \ldots,
n$.
$y_i$, $i=1, \ldots, n$, are generated through model (\ref
{model}). We are
going to apply the proposed estimation method based on the generated
$(s_i, X_i^{\mathrm{T}}, y_i)$, $i=1, \ldots, n$, to estimate
$\beta
_1(\cdot)$,
$\beta_2(\cdot)$, $\beta_3(\cdot)$, $\alpha$ and $\sigma^2$, and
examine the
accuracy of the proposed estimation procedure.

We use the Epanechnikov kernel $K(t) = 0.75 (1-t^2)_+$ as the kernel function
in the estimation procedure. The bandwidth used in the estimation is $0.4$.

We use mean squared error (MSE) to assess the accuracy of an estimator
of an
unknown constant parameter, mean integrated squared error (MISE) to
assess the
accuracy of an estimator of an unknown function.

\begin{table}[t]
\tabcolsep=0pt
\caption{The MISEs and MSEs}\label{tab1}
\begin{tabular*}{\tablewidth}{@{\extracolsep{\fill}}@{}lccccc@{}}
\hline
& $\bolds{\hat{\beta}_1(\cdot)}$ & $\bolds{\hat{\beta}_2(\cdot)}$ & $\bolds{\hat{\beta}_3(\cdot)}$ & $\bolds{\hat{\alpha}}$ & $\bolds{\hat{\sigma}^2}$\\
\hline
$n=400$ & 0.0769& 0.0642 & 0.0618 & 0.0128 & 0.0086 \\
$n=500$ & 0.0712& 0.0573 & 0.0539 & 0.0093 & 0.0065\\
$n=600$ & 0.0679& 0.0498 & 0.0474 & 0.0076 & 0.0053\\
\hline
\end{tabular*}
\tabnotetext[]{}{The column corresponding to the estimator of an unknown function
is the MISEs of the estimator for $n=400$, $n=500$ and $n=600$, corresponding
to the estimator of an unknown constant is the MSEs of the estimator.}
\end{table}

\begin{figure}[b]

\includegraphics{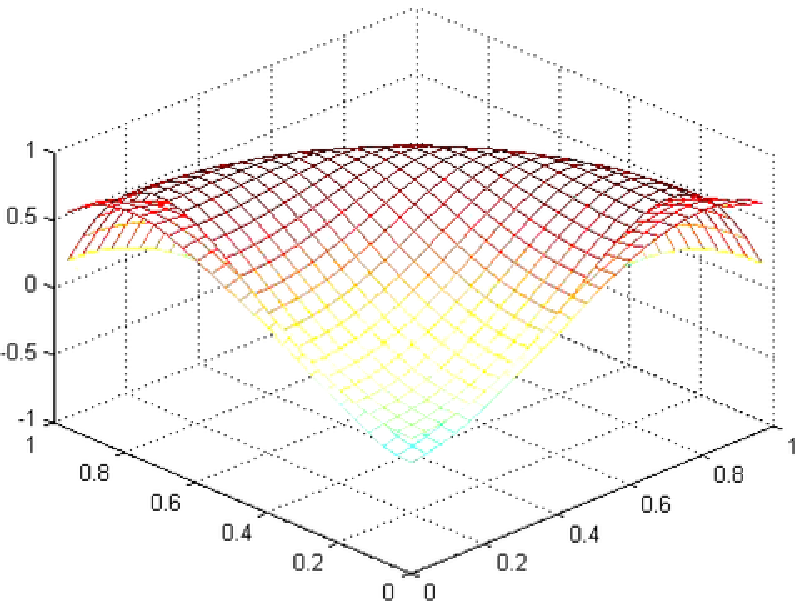}

\caption{The estimated $\beta_1(s)$ superimposed with $\beta _1(s)$.}\label{fig1}
\end{figure}

\begin{figure}

\includegraphics{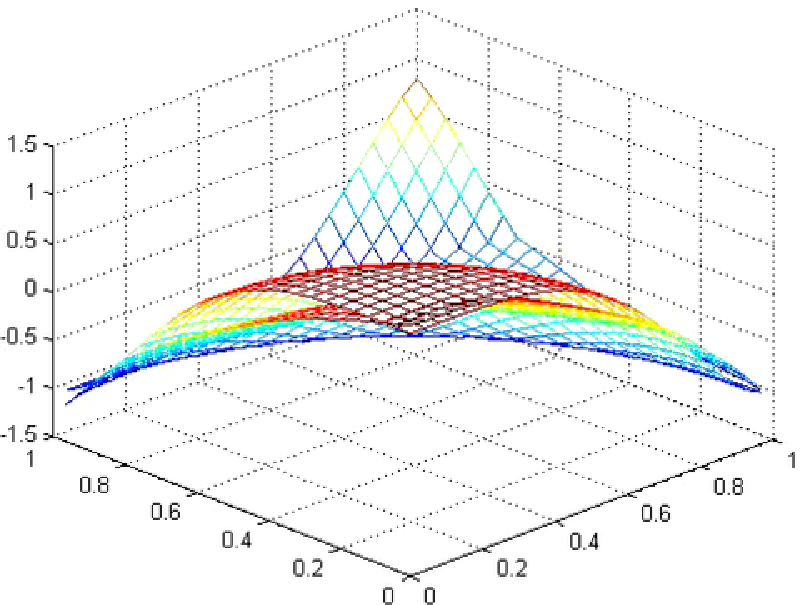}

\caption{The estimated $\beta_2(s)$ superimposed with $\beta_2(s)$.}\label{fig2}
\end{figure}

\begin{figure}

\includegraphics{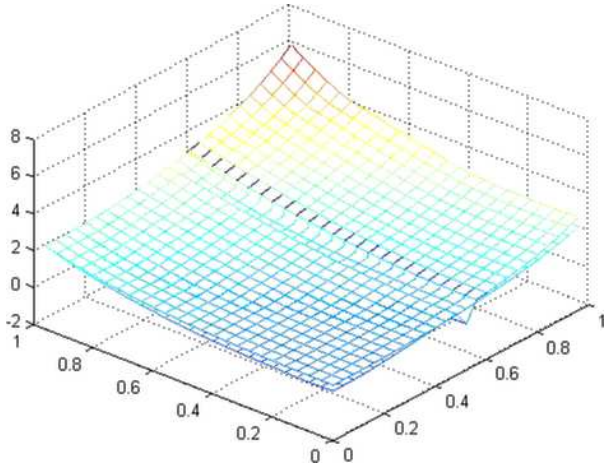}

\caption{The estimated $\beta_3(s)$ superimposed with $\beta_3(s)$.}\label{fig3}
\end{figure}

For each given sample size $n$, we do $200$ simulations. We compute the MSEs
of the estimators of the unknown constants and the MISEs of the
estimators of
the unknown functions for sample size $n=400$, $n=500$ and $n=600$. The
obtained
results are presented in Table~\ref{tab1}. Table~\ref{tab1} shows the proposed
estimation procedure works very well. To have a more visible idea about the
performance of the proposed estimation procedure, we set sample size $n=500$
and do $200$ simulations. We single out the one with median performance among
the $200$ simulations. The estimate of $\alpha$ coming from this simulation
is $ 0.407 $, the estimate of $\sigma^2$ is $ 0.976 $. The estimated
unknown functions from this simulation are presented in Figures~\ref{fig1},
\ref{fig2} and~\ref{fig3}, and are superimposed with the true functions.
All these show our estimation procedure works very well.
\end{exa}

\subsection{Performance of the model selection procedure}\label{sec5.2}

\begin{exa}\label{exa2}
In model~(\ref{model}), we set $p=5$, $\beta_1(\cdot)$,
$\beta_2(\cdot)$ and $\beta_3(\cdot)$ the same as that in Example~\ref{exa1},
$\beta_4(\cdot)= \sin^2(\|s\|^2 \pi)$, $\beta_5(\cdot) = \beta_5
= 1$.
We generate $X_i$, $s_i$, $\varepsilon_i$, $y_i$
$i=1, \ldots, n$, in the same way as that in Example~\ref{exa1}, except
that $X_i$
is from\vadjust{\goodbreak} $N(\mathbf{0}_5, I_5)$. Based on the generated data, we are
going to apply
the proposed $\mathrm{AIC}$ or $\mathrm{BIC}$ to select the correct
model, and
examine the
performances of the proposed $\mathrm{AIC}$, $\mathrm{BIC}$ and the two
algorithms in
identifying the constant components in model~(\ref{model}).

We still use the Epanechnikov kernel as the kernel function in the model
selection, however, the bandwidth used is 0.2 for AIC and 0.3 for BIC, which
is smaller than that for estimation. In general, the bandwidth used for
model selection should be
smaller than that for estimation. In fact, we have tried different
bandwidths, it turned out any bandwidth in a reasonable range such as
$[0.15, 0.3]$ for AIC, $[0.2, 0.35]$ for BIC would do the job very well.

\begin{table}[t]
\tabcolsep=0pt
\caption{Ratios of picking up each model in model selection}\label{tab2}
\begin{tabular*}{\tablewidth}{@{\extracolsep{\fill}}@{}lcccd{1.2}d{1.2}d{1.2}@{}}
\hline
& $\bolds{\{5\}}$ & $\bolds{\{1, 5\}}$ & $\bolds{\{4, 5\}}$
& \multicolumn{1}{c}{$\bolds{\{1, 4, 5\}}$} & \multicolumn{1}{c}{$\bolds{\{1, 2, 4, 5\}}$} & \multicolumn{1}{c@{}}{$\bolds{\{1, 2, 3, 4, 5\}}$}
\\
\hline
$n=400$ & 0.83 & 0.05& 0.07 & 0.02& 0.02& 0.01 \\
$n=500$ & 0.91 & 0.02& 0.05 & 0.02& 0& 0\\
$n=600$ & 0.94 & 0.02 & 0.03 & 0& 0 & 0.01\\
$n=400$ & 0.81 & 0.06& 0.08 & 0.05& 0& 0\\
$n=500$ & 0.89& 0.03& 0.06 & 0.02& 0& 0\\
$n=600$ & 0.92 & 0.01 & 0.04 & 0.01& 0.02 & 0\\
$n=400$ & 0.86 & 0.04& 0.05& 0.03& 0.01& 0.01\\
$n=500$ & 0.93 & 0.02& 0.03 & 0.01& 0.01& 0\\
$n=600$ & 0.96 & 0.01 & 0.02 & 0.01& 0 & 0\\
$n=400$ & 0.84& 0.05& 0.07 & 0.02& 0.01& 0.01\\
$n=500$ & 0.88& 0.05& 0.05 & 0.01& 0.01& 0\\
$n=600$ & 0.93& 0.03 & 0.03 & 0.01& 0 & 0\\
\hline
\end{tabular*}
\tabnotetext[]{}{The ratios of picking up each candidate model in $200$ simulations for
different sample sizes. $\{i_1, \ldots, i_k\}$ stands for the
model in which $\bolds{\beta}(\cdot)$ has its $i_1$th$, \ldots, i_k$th
components being
constant and the column corresponding to which is the ratios of picking up
this model among $200$ simulations. Row 2 to row 4 are the ratios obtained
based on AIC and backward elimination when sample size $n=400$, $n=500$ and
$n=600$. Row 5 to row 7 are the ratios obtained based on AIC and the CTAR
based algorithm, row 8 to row 10 are the ratios obtained based on BIC and
backward elimination, and row 11 to row 13 are the ratios obtained
based on BIC and the CTAR based algorithm.}
\end{table}

Due to the very expensive computation involved, for any given sample
size $n$,
we only do $200$ simulations, and in each
simulation, we apply either AIC or BIC coupled with either of the two proposed
algorithms to select model. For each candidate model, the ratios of
picking up
this model in the $200$ simulations are computed for different cases. The
results are presented
in Table~\ref{tab2}. We can see, from Table~\ref{tab2}, the proposed
BIC with
backward elimination performs best, and the others are doing reasonably
well, also.
\end{exa}

\section{Real data analysis}\label{sec6}
In this section, we are going to apply the proposed model~(\ref{model})
together with the proposed model selection and estimation method to analyze
the Boston house price data. Specifically, we are going to explore how some
factors such as the per capita crime rate by town (denoted by CRIM), average
number of rooms per dwelling (denoted by RM), index of accessibility to radial
highways (denoted by RAD), full-value property-tax rate per \$10,000 dollar
(denoted by TAX) and the percentage of the lower status of the population
(denoted by LSTAT) affect the median value of owner-occupied homes in
\$1000's (denoted by MEDV), and whether the effects of these factors
vary over
location or not.

\begin{table}[b]
\tabcolsep=0pt
\tablewidth=190pt
\caption{Estimates of the unknown constant coefficients}\label{tab5}
\begin{tabular*}{\tablewidth}{@{\extracolsep{\fill}}@{}lcc@{}}
\hline
$\bolds{\hat{\alpha}}$ & $\bolds{\hat{\beta}_3}$ & $\bolds{\hat{\beta}_5}$
\\
\hline
0.2210 & 0.3589 & $-0.4473$\\
\hline
\end{tabular*}
\end{table}

We use model~(\ref{model}) to fit the data with $y_i$, $x_{i1}$, $x_{i2}$,
$x_{i3}$, $x_{i4}$ and $x_{i5}$ being MEDV, CRIM, RM, RAD, TAX and LSTAT,
respectively, and $X_i = (x_{i1}, \ldots, x_{i5})^{\mathrm{T}}$.
The kernel
function used in either estimation procedure or model selection is
taken to
be the Epanechnikov kernel.

We first try to find which factors have location varying effects on the house
price, and which factors do not. This is equivalent to identifying the
constant coefficients in the model used to fit the data. We apply the proposed
BIC coupled with backward elimination to do the model selection, and the
bandwidth used is chosen to be $17\%$ of the range of the locations. The
obtained result shows the coefficients of $x_{i3}$ and $x_{i5}$ are constant,
which means all factors, except RAD and LSTAT, have location varying
effects on
the house price.

\begin{figure}

\includegraphics{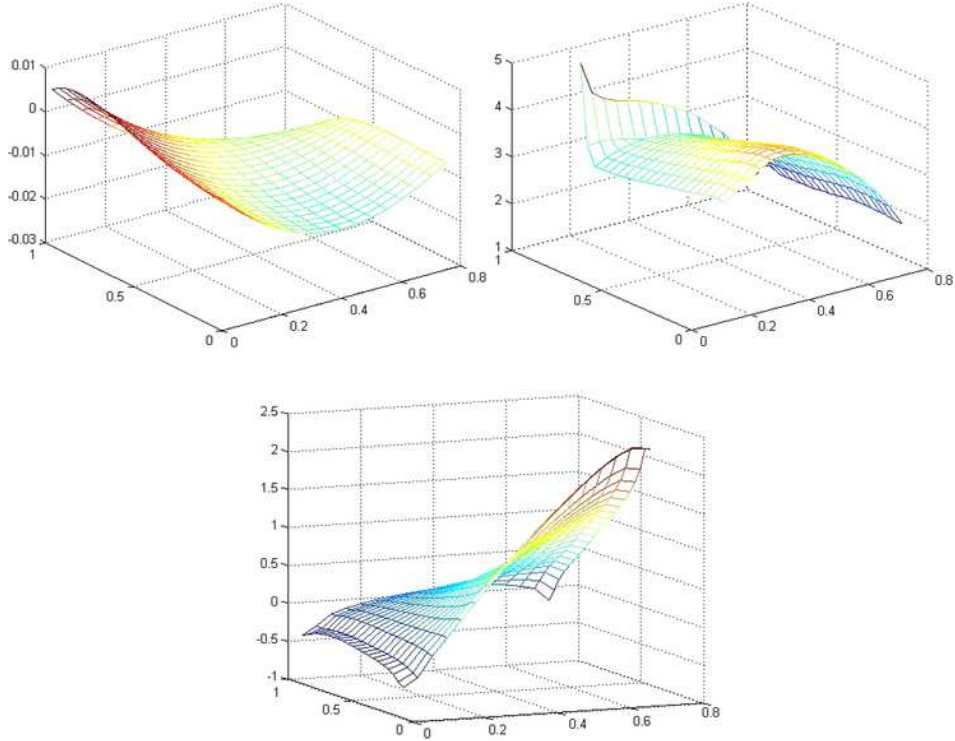}

\caption{The 3D plots of $\hat{\beta}_1(s)$, $\hat{\beta}_2(s)$
and $\hat{\beta}_4(s)$. The left one in the upper panel is
$\hat{\beta}_1(s)$, right one in the upper panel is $\hat{\beta}_2(s)$, and
the one in the lower panel is $\hat{\beta}_4(s)$.}\label{fig4}\vspace*{6pt}
\end{figure}

We now apply the chosen model
%
\begin{eqnarray}\label{real}
y_i &= & \alpha\sum_{j \neq i}
w_{ij} y_j + x_{i1} \beta_1(s_i)
\nonumber\\[-9pt]\\[-9pt]
&&{} + x_{i2} \beta_2(s_i) + x_{i3}
\beta_3 + x_{i4} \beta_4(s_i) +
x_{i5} \beta_5 + \varepsilon_i,\nonumber
\end{eqnarray}
$i = 1, \ldots, n$, where $w_{ij}$ is defined by~(\ref{weight}),
to fit
the data. The sample size of this data set is $n = 506$. The proposed
estimation procedure is used to estimate the unknown
functions and constants, and the bandwidth used in the estimation
procedure is
taken to be $60\%$ of the range of the locations. The estimates of the
unknown constants are presented in Table~\ref{tab5}, and the estimates
of the
unknown functions are presented in Figure~\ref{fig4}.

\begin{figure}

\includegraphics{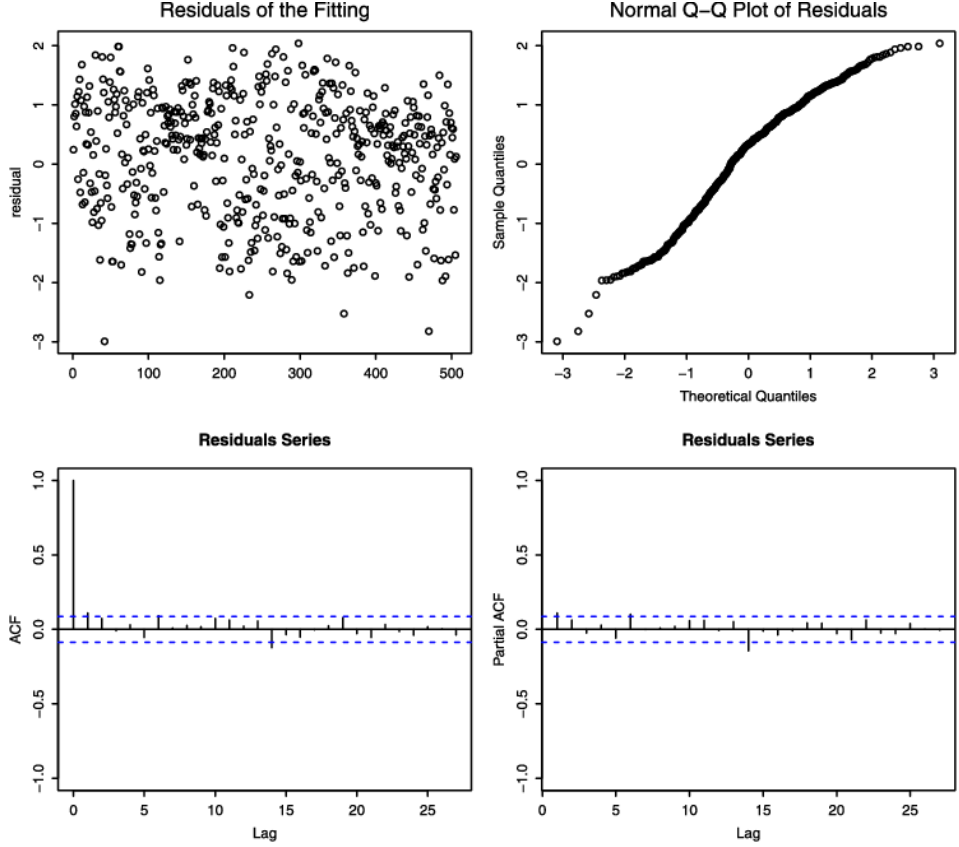}

\caption{The plot, normal $Q$--$Q$ plot, ACF and partial ACF of the
residuals of the fitting of (\protect\ref{real}) to the Boston house
price data.}\label{ex1}
\end{figure}

To see how well model~(\ref{real}) fits the data, we conduct some residual
analysis. The plot, normal $Q$--$Q$ plot, ACF and partial ACF of the
residuals of
the fitting are presented in Figure~\ref{ex1}. Figure~\ref{ex1} shows model
(\ref{real}) fits the data well.

As $\beta_3$ and $\beta_5$ can be interpreted as the impacts of RAD
and LSTAT,
respectively, Table~\ref{tab5} shows the index of accessibility to radial
highways has positive impact on house price and the percentage of the lower
status of the population has negative impact on house price.
Apparently, this
makes sense. Table~\ref{tab5} also shows that the estimate of $\alpha
$ is
$0.221$, which is an unignorable effect, and indicates the house prices
in a
neighborhood do affect each other. This is a true phenomenon in real world.

From Figure~\ref{fig4}, we can see the impact $\beta_1(\cdot)$ of
the per capita
crime rate by town on house price is negative and is clearly varying over
location. The impact $\beta_2(\cdot)$ of the average number of rooms
per dwelling on house price is positive and is also varying over location.
It is interesting to see that the impact of the average number of rooms per
dwelling is lower in the area where the impact of crime rate is high than
the area where the impact of crime rate is low. This implies that the
crime rate is a dominate factor on the house price in the area where the
impact of crime rate is high. Figure~\ref{fig4} also shows the association
between the house price and the full-value property-tax rate is varying over
location, and it is generally positive, however, there are some areas where
this association is negative. We can also see that the impact of the average
number of rooms per dwelling is lower in the area, where the association
between the house price and the full-value property tax rate is strong, than
the area where the association is weak.

\setcounter{equation}{0}
\begin{appendix}\label{app}
\section*{Appendix: Conditions and sketch of theoretical proofs}
To avoid confusion of notation, we use $\alpha_0$ to denote the true
value of
$\alpha$ in this section. Further, we rewrite $A=I_n-\alpha W$ as
$A(\alpha)$
to emphasis its dependence on $\alpha$ and abbreviate $A(\alpha_0)$
as $A$.

The following regularity conditions are needed to establish the asymptotic
properties of the estimators.

\subsection*{Conditions}
\begin{longlist}[($\tilde{7}$)]
\item[(1)] The kernel function $K(\cdot)$ is a bounded positive,
symmetric and
Lipschitz continuous function with a compact support on R.
$h\rightarrow0$.

\item[(2)] $\{\beta_i(\cdot), i=1, \ldots, p\}$ have
continuous second
partial derivatives.

\item[(3)] $\{X_1, \ldots, X_n\}$ is an i.i.d. random\vspace*{1pt} sample and is
independent of $\{\varepsilon_1, \ldots,\break  \varepsilon_n\}$. Moreover,
$E(X_1 X_1^{\mathrm{T}})$ is positive definite, $E\|X_1\|^{2q}<\infty
$ and
$E|\varepsilon_1|^{2q}<\infty$ for some $q>2$.

\item[(4)] $\{s_i\}$ is a sequence of fixed design points on a bounded compact
support $\mathcal{S}$. Further, there exists a positive joint density function
$f(\cdot)$ satisfying a Lipschitz condition such that
\[
\sup_{s\in\mathcal{S}} \Biggl|\frac{1}{n }\sum
_{i=1}^n \bigl[r(s_i)K_h \bigl(\|s_i-s\| \bigr)\bigr] - \int r(t)
K_h \bigl(\|t-s\| \bigr)f(t) \,dt \Biggr|=O(h)
\]
for any bounded continuous function $r(\cdot)$ and $K_h(\cdot
)=K(\cdot/h)/h^2$
where $K(\cdot)$ satisfies condition~(1). $f(\cdot)$ is bounded away
from zero
on $\mathcal{S}$.

\item[(5)] $w_{ii}=0$ for any $i$, and there exists a sequence $\rho
_n>0$ such
that $w_{ij} = O(1/\rho_n)$ uniformly with respect to $i$ and $j$.
Furthermore, the matrices $W$ and $A^{-1}$ are uniformly bounded in
both row
and column sums.

\item[(6)] $A^{-1}(\alpha)$ are uniformly bounded in either row or column
sums, uniformly in $\alpha$ in a compact support $\Delta$. The true
$\alpha_0$
is an interior point in~$\Delta$.

\item[(7)] $\lim_{n\rightarrow\infty} \frac{1}{n}
E[(G\mathbf{m}- Z)^{\mathrm{T}}(G\mathbf{m}-Z)]=\lambda_1 >0$.

\item[($\tilde{7}$)] $\lambda_1 =0$.

\item[(8)] $\rho_n$ is bounded and for any $\alpha\neq\alpha_0$,
\[
\lim_{n\rightarrow\infty} \biggl\{ \frac{1}{n}\log\bigl|\sigma^2
A^{-1} \bigl(A^{-1} \bigr)^{\mathrm{T}} \bigr| -
\frac{1}{n} \log\bigl|\sigma_a^2(\alpha)A^{-1}(
\alpha) \bigl(A^{-1}(\alpha) \bigr)^{\mathrm{T}} \bigr| \biggr\}\neq0,
\]
where $\sigma_a^2(\alpha)=\frac{\sigma^2}{n}
\operatorname{tr}\{(A(\alpha)A^{-1})^{\mathrm{T}} A(\alpha)A^{-1}\}$.

\item[(9)] $\rho_n \rightarrow\infty$, the row sums of $G$ have the uniform
order $O(1/\sqrt{\rho_n})$ and
\[
\lim_{n\rightarrow\infty} \frac{\rho_n}{n} E \bigl[(G\mathbf{m}-SG
\mathbf{m})^{\mathrm{T}}(G\mathbf{m}-SG\mathbf{m}) \bigr]= \lambda_4>0.
\]
\end{longlist}

\begin{rem}\label{rem1} Conditions (1)--(3) are commonly seen in nonparametric
estimation. They are not the weakest possible ones, but they are imposed
to facilitate the technical proofs. Since the sampling units can be regarded
as given, the fixed bounded design condition~(4) is made for technical
convenience. Of course, as in Linton \cite{r12}, condition~(4) does not
preclude $\{s_i\}_{i=1}^n$ from being generated by some random
mechanism. For
example, if $s_i$'s were i.i.d. with joint density $f(\cdot)$, then
condition~(4) holds with probability one which can be obtained in a
similar way
to Hansen
\cite{r7}. So, we can obtain our results by firstly conditional on
$\{s_i\}_{i=1}^n$, then some standard arguments.
\end{rem}

\begin{rem}\label{rem2}
Conditions (5)--(8) parallel the corresponding
conditions of Lee \cite{r10} and Su and Jin \cite{r14}. Conditions (5)--(6)
concern the essential features of the weight matrix for the model.
Condition~(7) is a sufficient condition which ensures that the likelihood
function of $\alpha$ has a unique maximizer. When
condition~($\tilde{7}$)\vadjust{\goodbreak}
holds and the elements of $W$ are uniformly bounded, the uniqueness of the
maximizer can be guaranteed by condition~(8). These two kinds of conditions
ensure that $\Omega$ which is the limit of the information matrix of
the finite-dimensional parameters is nonsingular. So, they are the crucial
conditions for
$\sqrt{n}$-rate of convergence of the finite-dimensional parameter estimators.
\end{rem}

\begin{rem}\label{rem3}
When $\rho_n\rightarrow\infty$, $\Omega$ is
nonsingular only when condition~(7) holds. Under condition~($\tilde{7}$),
$\Omega$ will become singular. The singularity of the matrix may have
implications on the rate of convergence of the estimators.
Nevertheless, we
follow Lee \cite{r10} and Su and Jin \cite{r14} to consider the
situation where
\[
\lim_{n\rightarrow\infty} \frac{\rho_n}{n} E \bigl[(G\mathbf{m}
)^{\mathrm{T}}(I_n -S)^{\mathrm{T}}(I_n-S) G
\mathbf{m} \bigr] = \lambda_4 \in(0, \infty).
\]
In this case, it is natural to assume that the elements of
$(I_n-S)G\mathbf{m}$ have the uniform order $O_P(1/\sqrt{\rho_n})$
which can be satisfied
by the assumption that the row sums of
$G$
are of uniform order $O(1/\sqrt{\rho_n})$.
\end{rem}

In the following, let $H$ be a diagonal matrix of size $3p$ with its
first $p$
elements on the diagonal\vspace*{1.5pt} being $1$ and the remaining elements being $h$,
$P=(I_n-S)^{\mathrm{T}}(I_n-S)$ and $\bolds{\varepsilon}=(\varepsilon
_1,\ldots,\varepsilon
_n)^{\mathrm{T}}$.
Moreover, like $\alpha_0$, we use $\sigma_0^2$ to denote the true
value of
$\sigma^2$ to avoid confusion of notation. Since the following notations
will be frequently used in the proofs, we list here for easy
reference:
\begin{eqnarray*}
l\bigl(\alpha,\sigma^2\bigr)&=& - \frac{n}{2} \log\bigl(\sigma^2\bigr)
+ \log\bigl(\bigl|A(\alpha)\bigr|\bigr)-\frac{1}{2\sigma^2}\bigl(A(\alpha) Y\bigr)^{\mathrm{T}} P
A(\alpha) Y,
\\
l_c(\alpha)&=&- \frac{n}{2}\log\tilde{\sigma}^2(\alpha) + \log\bigl|A(\alpha)\bigr|,
\\
\tilde{\sigma}^2(\alpha)&=&\frac{1}{n}\bigl(A(\alpha)Y\bigr)^{\mathrm{T}}P A(\alpha)Y,
\\
\bar{\sigma}^2(\alpha)&=&\frac{1}{n}E\bigl[\bigl(A(\alpha)Y\bigr)^{\mathrm{T}}P A(\alpha)Y\bigr],
\\
\sigma_a^2(\alpha)&=&\frac{\sigma_0^2}{n}\operatorname{tr}\bigl\{\bigl(A(\alpha)A^{-1}\bigr)^{\mathrm{T}} A(\alpha)A^{-1}\bigr\}.
\end{eqnarray*}

To prove the theorems,
the following lemmas are needed. Their
proofs and the more detailed proofs of the theorems
can be found in the supplementary material (Sun et~al. \cite{r21}).

\begin{lem}\label{lem1}
Let $\{ Y_i\} $
be a sequence of independent random variables and $\{s_i\}\in R^2$ are nonrandom
vectors. Suppose that for
some $q>2$, $\max_{i} E|Y_i|^{q}<\infty$.
Then under condition~(1), we have
\[
\sup_{s\in\mathcal{S}} \Biggl| \frac{1}{n}\sum
_{i=1}^n \bigl[ K_h\bigl(\|s_i-s\|\bigr) Y_i -E \bigl\{K_{h}\bigl(\|s_i-s\|\bigr) Y_i \bigr\} \bigr] \Biggr| = O_p
\biggl( \biggl\{\frac{\log n}{nh^2} \biggr\}^{1/2} \biggr),
\]
provided that $ n^{1-2/q} h^2/\log^2 n \rightarrow\infty$
and $\lim_{n\rightarrow\infty}\frac{1}{n}
\sum_{i=1}^n K_h(\|s_i-s\|)<\infty$
for any $s\in\mathcal{S}$.
\end{lem}

\begin{lem}\label{lem2}
Under conditions (1)--(4), when
$n^{1/2}h^2/\log^2 n\rightarrow\infty$,
\begin{longlist}[(1)]
\item[(1)] \mbox{}
\begin{eqnarray*}
\\[-40pt]
&& n^{-1}H^{-1}\mathcal{X}^{\mathrm{T}}\mathcal{W}\mathcal{X}H^{-1}
\\
&&\qquad = \pmatrix{
\kappa_0 f(s)\Psi& \mathbf{0}_{p\times2p} \vspace*{6pt}\cr
\mathbf{0}_{2p\times p} & \kappa_2 f(s)\Psi\otimes I_2}
+ O_P\bigl( c_n \mathbf{1}_{3p}\mathbf{1}^{\mathrm{T}}_{3p}\bigr)
\end{eqnarray*}
holds uniformly in $s\in\mathcal{S}$ where $c_n =h+\{\frac{\log n}{n
h^2}\}^{1/2}$,\vspace*{6pt}

\item[(2)]
\begin{eqnarray*}
\\[-40pt]
&& \bolds{\beta}(s)-(I_p, \mathbf{0}_{p\times2p})\bigl(\mathcal
{X}^{\mathrm{T}}\mathcal{W}\mathcal{X}\bigr)^{-1}
\mathcal{X}^{\mathrm{T}}\mathcal{W}\mathbf{m}
\\
&&\qquad = -\frac{\kappa_2 h^2}{ 2\kappa_0}
\bigl\{\bolds{\beta}_{uu}(s) + \bolds{\beta}_{vv}(s) \bigr\} +
o_p \bigl( h^2 \mathbf{1}_p \bigr)
\end{eqnarray*}
holds uniformly in $s\in\mathcal{S}$.
\end{longlist}
\end{lem}

\begin{lem}\label{lem3}
Under conditions (1)--(5), when
$n^{1/2}h^2/\log^2 n\rightarrow\infty$,
\[
n^{-1}H^{-1}\mathcal{X}^{\mathrm{T}}\mathcal{W}G\mathbf{m}
-n^{-1} E \bigl( H^{-1}\mathcal{X}^{\mathrm{T}}\mathcal{W}G
\mathbf{m} \bigr) = o_P( 1)
\]
uniformly in $s\in\mathcal{S}$.
\end{lem}

\begin{lem}\label{lem4}
Under conditions (1), (3), (4) and (5), when
$n^{1/2}h^2/\log^2 n\rightarrow\infty$, we have
(1) $ \frac{1}{n}E[\operatorname{tr}(P)] = 1+o(1)$,
(2) $\frac{1}{n} E[\operatorname{tr}(G^{\mathrm{T}}P) -\operatorname
{tr}(G)]= o(1)$,
(3)~$ \frac{1}{n}E[\operatorname{tr}(G^{\mathrm{T}}PG)-\operatorname
{tr}(G^{\mathrm{T}}G)]= o(1)$.
Further, when $nh^2 / \rho_n \rightarrow\infty$,\break
(4) $ \frac{\rho_n}{n} E[\operatorname{tr}(P)-n ] = o(1)$,
(5) $ \frac{\rho_n}{n} E[\operatorname{tr}(G^{\mathrm
{T}}P)-\operatorname{tr}(G)]= o(1)$,\break
(6) $ \frac{\rho_n}{n}E[\operatorname{tr}(G^{\mathrm
{T}}PG)-\operatorname{tr}(G^{\mathrm{T}}G)]=o(1)$.
\end{lem}

\begin{lem}\label{lem5}
Under conditions (1)--(5),
when $n^{1/2}h^2/\log^2 n\rightarrow\infty$,\break
(1) $
(G\mathbf{m})^{\mathrm{T}}P\mathbf{m} = o_P(n h^2)$.
Moreover, under the assumption that the second partial derivative of
$\bolds{\beta}(s)$ is Lipschitz\vspace*{1pt} continuous, we have (2)
$
(G\mathbf{m})^{\mathrm{T}}P\mathbf{m} = O_P(n h^3 + \{n h^2 \log n\}^{1/2})$.
\end{lem}

\begin{lem}\label{lem6}
Under conditions (1)--(5), when
$n^{1/2}h^2/\log^2 n\rightarrow\infty$ and\break  $nh^8\rightarrow0$, we have
(1) $n^{-1/2}L^{\mathrm{T}} P\mathbf{m} =o_P(1)$ for $L=\mathbf{m},
\bolds{\varepsilon}$
and $G\bolds{\varepsilon}$,\break
(2) $n^{-1} L^{\mathrm{T}} P G\mathbf{m} =o_P(1)$ for $L= \mathbf
{m}, \bolds{\varepsilon}
$ and $G\bolds{\varepsilon}$.
\end{lem}

\begin{lem}\label{lem7}
Under conditions (1)--(5), when
$n^{1/2}h^2/\log^2 n\rightarrow\infty$, we have
(1) $\frac{1}{n} \{ (G\mathbf{m})^{\mathrm{T}} P G\mathbf{m}
- E[ (G\mathbf{m})^{\mathrm{T}} P G\mathbf{m} ] \}=o_P(1)$,
(2) $\frac{1}{n} E[ (G\mathbf{m})^{\mathrm{T}} P G\mathbf{m} ]=
\frac{1}{n} E[ (G\mathbf{m}- Z)^{\mathrm{T}} ( G\mathbf{m} -Z) ]+o(1)$.
\end{lem}

\begin{lem}\label{lem8}
Under conditions (1)--(5), when
$n^{1/2}h^2/\log^2 n\rightarrow\infty$, we have
(1)~$ n^{-1/2}\{\bolds{\varepsilon}^{\mathrm{T}}P \bolds{\varepsilon}-
\bolds{\varepsilon}^{\mathrm{T}}\bolds{\varepsilon}\}=o_P(1) $,
(2)
$ n^{-1/2}\{\bolds{\varepsilon}^{\mathrm{T}}G^{\mathrm{T}} P \bolds
{\varepsilon}
- \bolds{\varepsilon}^{\mathrm{T}}G^{\mathrm{T}}\bolds{\varepsilon}\}
=o_P(1) $,
(3)
$ n^{-1/2}\{\bolds{\varepsilon}^{\mathrm{T}}G^{\mathrm{T}} P G \bolds
{\varepsilon}
- \bolds{\varepsilon}^{\mathrm{T}}G^{\mathrm{T}} G \bolds{\varepsilon}\}
=o_P(1) $,
(4)
$ n^{-1/2}\{(G\mathbf{m})^{\mathrm{T}} P \bolds{\varepsilon}
- (G \mathbf{m} - SG\mathbf{m})^{\mathrm{T}}\bolds{\varepsilon}\}= o_P(1)$.
\end{lem}

\begin{lem}\label{lem9}
Suppose that $B=(b_{ij})_{1\le i,j \le n}$
is a sequence of symmetric matrices with row and column sums
uniformly bounded and its elements are also uniformly bounded. Let
$\sigma_{Q_n}^2 $ be the\vspace*{-2pt} variance of $Q_n$ where $Q_n= (G\mathbf{m}-
SG\mathbf{m})^{\mathrm{T}}\bolds{\varepsilon}+ \bolds{\varepsilon
}^{\mathrm{T}}B\bolds{\varepsilon}-\sigma_0^2 \operatorname{tr}(B)$. Assume\vspace*{-2pt}
that the
variance $\sigma_{Q_n}^2 $ is $O(n)$ with $\{\frac{\sigma_{Q_n}^2
}{n}\}$ bounded away from zero, then we have under conditions (1)--(5)
that $\frac{Q_n}{\sigma_{Q_n}}\stackrel{D}{\longrightarrow}N(0,1)$.
\end{lem}

\begin{lem}\label{lem10}
Under conditions (1)--(5), and the row sums
of matrix
$G$ having the uniform order $O(1/\sqrt{\rho_n})$ and
$n^{1/2}h^2/\log^2 n\rightarrow\infty$, we have\break
(1)~$(G\mathbf{m})^{\mathrm{T}}P\mathbf{m} =o_P(\rho_n^{-1/2} n h^2) $.
Moreover, if the second partial derivative of $\bolds{\beta}(s)$ is Lipschitz
continuous, then (2)
$
(G\mathbf{m})^{\mathrm{T}}P\mathbf{m}
= O_P( \rho_n^{-1/2} n h^3 +\break \{ n h^2 \log n/\rho_n \}^{1/2})$.
\end{lem}

\begin{lem}\label{lem11}
Under conditions (1)--(5) and the row sums
of matrix
$G$ having the uniform order $O(1/\sqrt{\rho_n})$, when
$n^{1/2}h^2/\log^2 n\rightarrow\infty$, $\rho_n \rightarrow\infty
$, $\rho_n
h^4\rightarrow0$ and $nh^2/\rho_n \rightarrow\infty$, we have (1)
$\frac{\rho_n}{n} \mathbf{m}^{\mathrm{T}} P \mathbf{m} =o_P(1)$,\break  (2)
$\frac{\rho_n}{n} L^{\mathrm{T}}
P G\mathbf{m} =o_P(1) $ for $L=\mathbf{m}, \bolds{\varepsilon}$ and
$G\bolds{\varepsilon}$, (3) $
\sqrt{\frac{\rho_n}{n} } (G\bolds{\varepsilon})^{\mathrm{T}} P
\mathbf{m} =o_P(1)$,\break  (4)
$\frac{\rho_n}{n} \{ (G\mathbf{m})^{\mathrm{T}} P G\mathbf{m}
- E[
(G\mathbf{m})^{\mathrm{T}} P G\mathbf{m}
] \}=o_P(1)$,
(5) $\sqrt{\frac{\rho_n}{n}}\{\bolds{\varepsilon}^{\mathrm
{T}}G^{\mathrm{T}} P \bolds{\varepsilon}
-\break  \bolds{\varepsilon}^{\mathrm{T}}G^{\mathrm{T}}\bolds{\varepsilon}\}
=o_P(1) $,
(6) $ \sqrt{\frac{\rho_n}{n}}\{\bolds{\varepsilon}^{\mathrm
{T}}G^{\mathrm{T}} P G \bolds{\varepsilon}
- \bolds{\varepsilon}^{\mathrm{T}}G^{\mathrm{T}} G \bolds{\varepsilon}\}
=o_P(1) $,\break
(7) $\sqrt{\frac{\rho_n}{n}}\{ (G\mathbf{m})^{\mathrm{T}} P \bolds
{\varepsilon}
- (G \mathbf{m} - SG\mathbf{m})^{\mathrm{T}}\bolds{\varepsilon}\}= o_P(1)$.
\end{lem}

\begin{lem}\label{lem12}
Suppose that $B=(b_{ij})_{1\le i,j \le n}$
is a sequence of symmetric matrices with row and column sums uniformly bounded.
Let $\sigma_{Q_n}^2 $ be the variance of $Q_n$ where
$Q_n= (G\mathbf{m}- SG\mathbf{m})^{\mathrm{T}}\bolds{\varepsilon}+
\bolds{\varepsilon}^{\mathrm{T}}B\bolds{\varepsilon}-\sigma
_0^2 \operatorname{tr}(B)$. Assume
that the variance $\sigma_{Q_n}^2 $ is $O(n/\rho_n)$ with
$\{\frac{\rho_n}{n}\sigma_{Q_n}^2 \}$ bounded\vspace*{1pt} away from zero, the elements
of $B$ are of uniform order $O(1/\rho_n)$ and the row sums of $G$ of uniform
order $O(1/\sqrt{\rho_n})$, we have under $\rho_n\rightarrow\infty
$ and
conditions (1)--(5) that $\frac{Q_n}{\sigma_{Q_n}}\stackrel
{D}{\longrightarrow}N(0,1)$.
\end{lem}

In the proofs of the theorems, we will use the facts that for constant
matrices $B=(b_{ij})$ and $D=(d_{ij})$,
$
\operatorname{var}(\bolds{\varepsilon}^{\mathrm{T}}B\bolds{\varepsilon
})=(\mu_4-3\sigma_0^4)\sum_{i=1}^n b_{ii}^2
+ \sigma_0^4 [\operatorname{tr}(B B^{\mathrm{T}}) +\operatorname{tr}(B^2)]
$
and
\[
E\bigl(\bolds{\varepsilon}^{\mathrm{T}}B\bolds{\varepsilon}\bolds{\varepsilon
}^{\mathrm{T}}D\bolds{\varepsilon}\bigr)=\bigl(\mu_4-3\sigma_0^4\bigr)\sum_{i=1}^n b_{ii}d_{ii}
+\sigma_0^4\bigl[\operatorname{tr}(B)\operatorname{tr}(D) + \operatorname
{tr}(B D) +\operatorname{tr}\bigl(B D^{\mathrm{T}}\bigr)\bigr].
\]

Moreover, we will frequently use the following facts by condition~(5)
(see Lee \cite{r10}) without being clearly pointed out:
\begin{longlist}[(2)]
\item[(1)] the elements of $G=W A^{-1}$ are $O(1/\rho_n)$
uniformly with respect to $i$ and $j$.

\item[(2)] The matrix $G=W A^{-1}$
is uniformly bounded in both row and column sums.
\end{longlist}

\begin{pf*}{Proof of Theorem~\ref{teo1}}
We will first show that $\Omega$ is nonsingular.
Let $\mathbf{d}=(d_1,d_2)^{\mathrm{T}}$ be a constant vector such
that $\Omega\mathbf{d}
=\mathbf{0}_2$.
Then it is sufficient to show that $\mathbf{d}=\mathbf{0}_2$. From
the second
equation of
$\Omega\mathbf{d}=\mathbf{0}_2$, we have that $d_2=-2\sigma_0^2\lim_{n\rightarrow
\infty}\frac{1}{n}\operatorname{tr}(G)d_1$. Plugging $d_2$ into the
first equation of
$\Omega\mathbf{d}=\mathbf{0}_2$, we have that
\[
d_1 \biggl\{ \frac{1}{\sigma_0^2}\lambda_1 + \lim
_{n\rightarrow\infty} \biggl[ \frac{1}{n}\operatorname{tr} \bigl(
\bigl(G+G^{\mathrm{T}} \bigr)G \bigr)-\frac
{2}{n^2}\operatorname{tr}^2(G)
\biggr] \biggr\}=0.
\]
It follows by condition~(7) that $ \lambda_1 >0 $. Moreover,
$\operatorname{tr}\{(G+G^{\mathrm{T}})G\}-\frac{2}{n}\operatorname{tr}^2(G)=
\frac{1}{2}\operatorname{tr}\{(\widetilde{G}{}^{\mathrm{T}}+\widetilde{G})
(\widetilde{G}{}^{\mathrm{T}}+\widetilde{G})^{\mathrm{T}}\}\ge0$ where
$\widetilde{G}=G-\frac{1}{n}\operatorname{tr}(G)I_n$.
As we have by condition~(5) that
$\operatorname{tr}\{(\widetilde{G}{}^{\mathrm{T}}+
\widetilde{G})(\widetilde{G}{}^{\mathrm{T}}+\widetilde{G})^{\mathrm
{T}}\}
=O(\frac{n}{\rho_n})$,
if condition~($\tilde{7}$) holds, condition~(8) implies that the limit of
$\frac{1}{2n}\operatorname{tr}\{(\widetilde{G}{}^{\mathrm{T}}+
\widetilde{G})(\widetilde{G}{}^{\mathrm{T}}+\widetilde{G})^{\mathrm
{T}}\}
> 0$.
Therefore, $d_1=0$ and $d_2=0$.

Next, we will follow the idea of Lee \cite{r10} to show the
consistency of
$\hat{\alpha}$. Define $Q(\alpha)$ to be
$\max_{\sigma^2} E [l(\alpha,\sigma^2) ]$
by ignoring the constant term. The optimal solution of this maximization
problem is
$\bar{\sigma}^2(\alpha)
=\frac{1}{n}E[(A(\alpha)Y)^{\mathrm{T}}P A(\alpha)Y]$.
Consequently,
\[
Q(\alpha)=-n/2\cdot\log\bar{\sigma}^2(\alpha) + \log\bigl|A(\alpha)\bigr|.
\]

According\vspace*{1pt} to White (\cite{r18}, Theorem 3.4),
it suffices to show the
uniform convergence of
$n^{-1}\{l_c(\alpha)-Q(\alpha)\}$ to zero in
probability on $\Delta$ and
the unique maximizer condition
that
%
\begin{equation}\label{id}
\lim\sup_{n\rightarrow\infty}\max_{\alpha\in
N^c(\alpha_0,\delta)}n^{-1}
\bigl[Q( \alpha)-Q(\alpha_0) \bigr]<0\qquad\mbox{for any } \delta>0,
\end{equation}
where $N^c(\alpha_0,\delta)$ is the complement of an open
neighborhood of $\alpha_0$ in $\Delta$ with diameter $\delta$.

Note that
$\frac{1}{n}l_c(\alpha)-\frac{1}{n}Q(\alpha)=-\frac{1}{2}\{\log
\tilde{\sigma}^2(\alpha)-\log\bar{\sigma}^2(\alpha)\}$, then to
show the
uniform convergence, it is sufficient to show that
$\tilde{\sigma}^2(\alpha)-\bar{\sigma}^2(\alpha)=o_P(1)$ uniformly
on $\Delta$ and $\bar{\sigma}^2(\alpha)$ is uniformly
bounded away from zero
on $\Delta$.

As
$A(\alpha)A^{-1}=I_n + (\alpha_0 -\alpha)G$ by $W A^{-1}=G$,
the result $\tilde{\sigma}^2(\alpha)-\bar{\sigma}^2(\alpha
)=o_P(1)$ uniformly
on $\Delta$ can be obtained by
straightforward calculations, Lemmas~\ref{lem4}(1)--(3),
\ref{lem6}, \ref{lem7}(1),
\ref{lem8}(1)--(3) and Chebyshev's inequality.

Now we will show that $\bar{\sigma}^2(\alpha)$ is bounded away from
zero uniformly on $\Delta$. As we know by simple calculations and
Lemma~\ref{lem4}(1)--(3) that
%
\begin{equation}
\bar{\sigma}^2(\alpha) \ge\sigma_0^2
n^{-1} \operatorname{tr} \bigl\{ \bigl(A(\alpha)A^{-1}
\bigr)^{\mathrm{T}}A(\alpha)A^{-1} \bigr\}+o(1), \label{sigbar}
\end{equation}
it suffices to show that
$\sigma_a^2(\alpha)=\frac{\sigma_0^2}{n}
\operatorname{tr}\{(A(\alpha)A^{-1})^{\mathrm{T}}A(\alpha)A^{-1}\}$
is uniformly bounded away from zero on $\Delta$. To do so, we define an
auxiliary spatial autoregressive (SAR) process: $Y=\alpha_0WY+\bolds
{\varepsilon}$ with
$\bolds{\varepsilon}\sim N(\mathbf{0}, \sigma_0^2 I_{n})$. Its log
likelihood function
without the
constant term is
\[
l_a \bigl(\alpha,\sigma^2 \bigr)= - \frac{n}{2}
\log\sigma^2 +\log\bigl|A(\alpha)\bigr| - \frac{1}{2\sigma^2} \bigl(A(\alpha)Y
\bigr)^{\mathrm{T}}A(\alpha)Y.
\]
Set $Q_a(\alpha)$ to be $\max_{\sigma^2}E_a[l_a(\alpha,
\sigma^2)]$
by ignoring the constant term, where $E_a$ is the expectation under
this SAR
process. It can be easily shown that
\[
Q_a(\alpha) =-n/2\cdot\log\sigma_a^2(
\alpha) + \log\bigl|A(\alpha)\bigr|.
\]
So, we have by Jensen's inequality that $Q_a(\alpha)\le Q_a(\alpha_0)$
for all
$\alpha\in\Delta$, hence it follows:
\[
-\frac{1}{2}\log\sigma_a^2(\alpha)\le-
\frac{1}{2}\log\sigma_0^2 + \frac{1}{n}
\bigl(\log\bigl|A(\alpha_0)\bigr| - \log\bigl|A(\alpha)\bigr| \bigr)
\]
uniformly on $\Delta$. Since we have, by the mean value theorem
and conditions~\mbox{(5)--(6)}, that
$
n^{-1}\{\log|A(\alpha_2)|-\log|A(\alpha_1)|\}=O(1)
$
uniformly\vspace*{2pt} in $ \alpha_1$ and $\alpha_2$ on $\Delta$, it follows that
$-\frac{1}{2}\log\sigma_a^2(\alpha)$ is bounded from\vspace*{2pt} above for any
$\alpha\in\Delta$. Therefore, the statement that $\sigma_a^2(\alpha)$
is uniformly bounded away from zero on $\Delta$ can be established
by a counter argument.

To show the uniqueness condition~(\ref{id}), write
\begin{eqnarray*}
n^{-1} \bigl[Q(\alpha)-Q(\alpha_0) \bigr]& = &
n^{-1} \bigl[Q_a(\alpha)-Q_a(
\alpha_0) \bigr] + 2^{-1} \bigl[\log\sigma_a^2(
\alpha) - \log\bar{\sigma}^2(\alpha) \bigr]
\\
& &{} +2^{-1} \bigl[\log\bar{\sigma}^2(
\alpha_0) - \log\sigma_0^2 \bigr],
\end{eqnarray*}
it follows, by Lemmas~\ref{lem4}(1) and~\ref{lem6}(1) and $\bar{\sigma}^2(\alpha_0)$ being
bounded away from zero, that
$\log\bar{\sigma}^2(\alpha_0) - \log\sigma_0^2=o(1)$. Moreover,
we have
already shown in~(\ref{sigbar}) that
$\lim_{n\rightarrow\infty}
[\sigma_a^2(\alpha) - \bar{\sigma}^2(\alpha) ]\le0$,
hence,
\[
\lim\sup_{n\rightarrow\infty} \max_{\alpha\in N^c (\alpha_0, \delta)} n^{-1}
\bigl[Q( \alpha)-Q(\alpha_0) \bigr]\le0\qquad\mbox{for any }
\delta>0.
\]

Now we will show that the above inequality holds strictly. It can be shown
that $n^{-1}Q(\alpha)$ is uniformly equicontinuous in $\alpha$ on
$\Delta$ by
Lemmas~\ref{lem4}(1)--(3), \ref{lem6}~and~\ref{lem7}(2) and the mean value theory. By the
compactness of $N^c(\alpha_0, \delta)$, there exists an $\delta>0$
and a
sequence $\{\alpha_n\} $ in $N^c(\alpha_0, \delta)$ converging to a point
$\alpha^{*}\neq\alpha_0$ such that
$\lim_{n\rightarrow\infty}n^{-1}[Q(\alpha_n)-Q(\alpha_0)]= 0$.
Because
$\lim_{n\rightarrow\infty} n^{-1}[Q(\alpha_n)-Q(\alpha
^{*})]= 0$ as
$\alpha_n\rightarrow\alpha^{*}$, it follows that
%
\begin{equation}
\lim_{n\rightarrow\infty} n^{-1} \bigl[Q \bigl(
\alpha^{*} \bigr)-Q(\alpha_0) \bigr]= 0. \label{q0}
\end{equation}

Since $Q_a(\alpha^{*})-Q_a(\alpha_0)\le0$ and
$\lim_{n\rightarrow\infty}
[\sigma_a^2(\alpha^{*})-\bar{\sigma}^2(\alpha^{*})]\le0$,
(\ref{q0}) is possible only if (i)
$ \lim_{n\rightarrow\infty}
[\sigma_a^2(\alpha^{*})-\bar{\sigma}^2(\alpha^{*})]
=0$
and
(ii)
$\lim_{n\rightarrow\infty}n^{-1}[Q_a(\alpha^{*})-Q_a(\alpha
_0)]= 0$.
However, (i) is a contradiction when condition~(7) holds by Lemmas~\ref{lem4}(1)--(3), \ref{lem6}~and~\ref{lem7}(2). If condition~($\tilde{7}$) holds, the contradiction follows from
(ii) by condition~(8).

The consistency of $\hat{\sigma}^2$ can be obtained straightforwardly
by Lemmas~\ref{lem4}(1)--(3), \ref{lem6},~\ref{lem7},~\ref{lem8}(1)--(3),
Chebyshev's inequality and $\hat{\alpha
}\stackrel{P}{\longrightarrow}\alpha_0$.
\end{pf*}

\begin{pf*}{Proof of Theorem~\ref{teo2}}
Denoting $\bolds{\theta}=(\alpha,\sigma^2)^{\mathrm{T}}$ and
$\bolds{\theta}_0=(\alpha_0,\sigma_0^2)^{\mathrm{T}}$, we get by
Taylor's
expansion that
\[
0=\frac{\partial l(\hat{\bolds{\theta}})}{\partial\bolds{\theta
}} =\frac
{\partial
l(\bolds{\theta}_0)}{\partial\bolds{\theta}}+\frac{\partial^2
l(\tilde{\bolds{\theta}})}{\partial\bolds{\theta}\,\partial\bolds{\theta
}^{\mathrm{T}}}(\hat{\bolds{
\theta}}-\bolds{\theta}_0),
\]
where $\tilde{\bolds{\theta}}=(\tilde{\alpha}, \tilde{\sigma
}^2)^{\mathrm{T}}$
lies between $\hat{\bolds{\theta}}$ and $\bolds{\theta}_0$,
and thus converges to $\bolds{\theta}_0$ in probability by Theorem~\ref{teo1}.
The asymptotic
distribution of $\hat{\bolds{\theta}}$ can be obtained by showing that
$
-\frac{1}{n}\frac{\partial^2 l(\tilde{\bolds{\theta}})}{\partial
\bolds{\theta}\,\partial\bolds{\theta}^{\mathrm{T}}}\stackrel
{P}{\longrightarrow
}\Omega$ and
$
\frac{1}{\sqrt{n}}\frac{\partial l(\bolds{\theta}_0)}{\partial
\bolds{\theta}}
\stackrel{D}{\longrightarrow}N(\mathbf{0}, \Sigma+\Omega)$,
where $\Omega$ is a nonsingular\vspace*{1pt} matrix by Theorem~\ref{teo1}.

By straightforward calculations, it can be easily obtained that
%
\begin{eqnarray}
\label{2der} \frac{1}{n}\frac{\partial^2 l(\bolds{\theta})}{\partial
\alpha^2} &=& -\frac{1}{n}
\operatorname{tr} \bigl( \bigl[WA^{-1}(\alpha) \bigr]^2
\bigr) - \frac{1}{\sigma^2 n}(WY)^{\mathrm{T}} P WY,
\nonumber
\\
\frac{1}{n} \frac{\partial^2 l(\bolds{\theta})} {
\partial\sigma^2\,\partial\sigma^2} &=& \frac{1}{2\sigma^4} - \frac
{1}{\sigma^6 n}
\bigl(A(\alpha)Y \bigr)^{\mathrm{T}}PA(\alpha)Y,
\\
\frac{1}{n} \frac{\partial^2 l(\bolds{\theta})} {
\partial\alpha\,\partial\sigma^2} &=&-\frac{1}{\sigma^4 n} (WY)^{\mathrm{T}}
P A(\alpha)Y.
\nonumber
\end{eqnarray}
As $A(\tilde{\alpha})A^{-1}=I_n + (\alpha_0 -\tilde{\alpha})G$
by $G=W A^{-1}$, we have
$
\frac{1}{n}
\frac{\partial^2 l(\tilde{\bolds{\theta}})} {
\partial\bolds{\theta}\,\partial\bolds{\theta}^{\mathrm{T}}}
-\frac{1}{n}
\frac{\partial^2 l(\bolds{\theta}_0)} {
\partial\bolds{\theta}\,\partial\bolds{\theta}^{\mathrm{T}}}
=o_P(1)
$
by Lemmas~\ref{lem6},~\ref{lem7},~\ref{lem8}(1)--(3), Chebyshev's inequality, mean value theorem
and $\tilde{\bolds{\theta}}\stackrel{P}{\longrightarrow}\bolds
{\theta}_0$. Furthermore, we have, by
Lemmas~\ref{lem6},~\ref{lem7},~\ref{lem8}(1)--(3) and Chebyshev's inequality that
$-\frac{1}{n}\frac{\partial^2 l(\bolds{\theta}_0)}{\partial
\bolds{\theta}\,\partial\bolds{\theta}^{\mathrm{T}}}\stackrel
{P}{\longrightarrow}\Omega$.

In the following, we will establish the asymptotic distribution of
$\frac{1}{\sqrt{n}}\frac{\partial l(\bolds{\theta}_0)}{\partial
\bolds{\theta}}$.
It follows by Lemma~\ref{lem5}(2) that $\frac{1}{\sqrt{n}}(G\mathbf
{m})^{\mathrm{T}
}P \mathbf{m} = o_P(1)$
when $nh^6\rightarrow0$ and $h^2 \log n \rightarrow0$. So, we have, by
straightforward calculations, Lemmas~\ref{lem6}(1) and~\ref{lem8}, that
\[
\frac{1}{\sqrt{n}}\frac{\partial l(\bolds{\theta}_0)}{\partial
\alpha} = \frac{1}{\sigma_0^2\sqrt{n}} \bigl[ (G \mathbf{m} -
SG\mathbf{m})^{\mathrm{T}}\bolds{\varepsilon}+ \bigl\{\bolds{\varepsilon
}^{\mathrm{T}}G \bolds{\varepsilon}-\sigma_0^2
\operatorname{tr}(G) \bigr\} \bigr]+o_P(1)
\]
and
\[
\frac{1}{\sqrt{n}}\frac{\partial l(\bolds{\theta}_0)}{\partial
\sigma^2} = \frac{1}{2\sigma_0^4\sqrt{n}} \bigl\{\bolds{
\varepsilon}^{\mathrm
{T}}\bolds{\varepsilon}-n\sigma_0^2
\bigr\}+o_P(1).
\]
By straightforward calculations and Lemma~\ref{lem7}(2), we have
$
E (\frac{1}{n}\frac{\partial
l(\bolds{\theta}_0)}{\partial\bolds{\theta}}
\frac{\partial
l(\bolds{\theta}_0)}{\partial\bolds{\theta}^{\mathrm{T}}} )
=\Sigma+ \Omega+o(1)$.

Finally,\vspace*{-2pt} as the components of $\frac{1}{\sqrt{n}}\frac{\partial
l(\bolds{\theta}_0)}{\partial\bolds{\theta}}= (
\frac{1}{\sqrt{n}}\frac{\partial l(\bolds{\theta}_0)}{\partial
\alpha},
\frac{1}{\sqrt{n}}\frac{\partial l(\bolds{\theta}_0)}{\partial
\sigma^2} )^{\mathrm{T}}$ are linear-quadratic forms of double
arrays, using
Lemma~\ref{lem9} we have
$
\frac{1}{\sqrt{n}}\frac{\partial l(\bolds{\theta}_0)}{\partial
\bolds{\theta}}
\stackrel{D}{\longrightarrow} N(\mathbf{0}, \Sigma+ \Omega)$.
\end{pf*}

\begin{pf*}{Proof of Theorem~\ref{teo3}}
It can be easily shown that
\begin{eqnarray*}
& & \sqrt{nh_1^2f(s)} \bigl(\hat{\bolds{
\beta}}(s)-\bolds{\beta}(s) \bigr)
\\
&&\qquad = \sqrt{nh_1^2f(s)}(I_p,
\mathbf{0}_{p\times2p}) \bigl(\mathcal{X}_1^{\mathrm{T}}
\mathcal {W}_1\mathcal{X}_1 \bigr)^{-1}
\mathcal{X}_1^{\mathrm{T}
}\mathcal{W}_1\bolds{
\varepsilon}
\\
&&\quad\qquad{} + \sqrt{nh_1^2f(s)}(
\alpha_0- \hat{\alpha}) (I_p, \mathbf{0}_{p\times2p})
\bigl( \mathcal{X}_1^{\mathrm{T}}\mathcal{W}_1
\mathcal{X}_1 \bigr)^{-1}\mathcal{X}_1^{\mathrm{T}}
\mathcal{W}_1 WY
\\
&&\quad\qquad{} + \sqrt{nh_1^2f(s)} \bigl
\{(I_p, \mathbf{0}_{p\times2p}) \bigl(\mathcal{X}_1^{\mathrm{T}}
\mathcal{W}_1\mathcal{X}_1 \bigr)^{-1}
\mathcal {X}_1^{\mathrm{T}}\mathcal{W}_1\mathbf{m} -
\bolds{ \beta}(s) \bigr\}
\\
&&\qquad \equiv J_{n1} + J_{n2} + J_{n3},
\end{eqnarray*}
where $\mathcal{X}_1$ and $\mathcal{W}_1$ are $\mathcal{X}$ and
$\mathcal{W}$ with $h$
being replaced
by $h_1$.

Let $H_1$ be $H$ with $h$ being replaced by $h_1$. It follows by
straightforward calculations that
\begin{eqnarray*}
& & n^{-1}h_1^2f(s)\operatorname{cov} \bigl
\{H_1^{-1}\mathcal{X}_1^{\mathrm
{T}}
\mathcal{W}_1\bolds{\varepsilon} \bigr\}
\\
&&\qquad  = \sigma_0^2
n^{-1}h_1^2f(s)E \bigl\{H_1^{-1}
\mathcal{X}_1^{\mathrm
{T}}\mathcal{W}_1^2
\mathcal{X}_1 H_1^{-1} \bigr\}
\\
&&\qquad = \sigma_0^2 f^2(s)\pmatrix{
\nu_0\Psi+o_P \bigl(\mathbf{1}_{p}
\mathbf{1}_p^{\mathrm{T}} \bigr) & o_P \bigl(
\mathbf{1}_{p}\mathbf{1}_{2p}^{\mathrm{T}} \bigr)
\vspace*{5pt}\cr
o_P \bigl(\mathbf{1}_{2p}\mathbf{1}_p^{\mathrm{T}}
\bigr) & \nu_2\Psi\otimes I_2 +o_P \bigl(
\mathbf{1}_{2p}\mathbf{1}_{2p}^{\mathrm{T}} \bigr)}
\end{eqnarray*}
this together with the central limit theorem, Lemma~\ref{lem2}(1) and Slutsky's theorem
lead to
\[
J_{n1}\stackrel{D} {\longrightarrow}N \bigl(\mathbf{0},
\nu_0\kappa_0^{-2}\sigma_0^2
\Psi^{-1} \bigr).
\]
It follows immediately from Lemmas~\ref{lem3},~\ref{lem2}(1) and condition~(4) that
\[
(I_p, \mathbf{0}_{p\times2p}) \bigl(\mathcal{X}_1^{\mathrm{T}}
\mathcal{W}_1\mathcal{X}_1 \bigr)^{-1}
\mathcal {X}_1^{\mathrm{T}}\mathcal{W}_1 G(\mathbf{m}
+ \bolds{\varepsilon})=O_P(1).
\]
When $nh_1^6=O(1)$ and $h/h_1\rightarrow0$, we have
$\sqrt{\frac{h_1^2}{n}}(G\mathbf{m})^{\mathrm{T}}P\mathbf{m} =o_P(1)$
using Lemma~\ref{lem5}(1). It can\vspace*{1pt}
be seen in the proof of Theorem~\ref{teo2} that
$\sqrt{nh_1^2}(\hat{\alpha}-\alpha_0)=o_P(1)$ under the assumptions of
Theorem~\ref{teo3}. Therefore, $J_{n2} =o_P(1)$.

The results of $J_{n1}$ and $J_{n2}$ together with Lemma~\ref{lem2}(2), $nh_1^6=O(1)$
and \mbox{$h/h_1\rightarrow0$} lead to the theorem.
\end{pf*}

\begin{pf*}{Proof of Theorem~\ref{teo4}}
It is obvious from the proof of nonsingularity of $\Omega$ in Theorem
\ref{teo1} that
$\Omega$ is singular under condition~(9).

Like Lee \cite{r10}, to prove the consistency of $\hat{\alpha}$, it suffices
to show that
\[
\frac{\rho_n}{n} \bigl\{ l_c(\alpha)- l_c(
\alpha_0)- \bigl[ Q(\alpha)- Q(\alpha_0) \bigr] \bigr
\}=o_P(1)\qquad\mbox{uniformly on } \Delta,
\]
where $Q(\alpha)=- n/2 \cdot\log\bar{\sigma}^2(\alpha)
+ \log|A(\alpha)|$ and $\alpha_0$ is the unique maximizer.

It follows by the mean value theorem that
\begin{eqnarray*}
&& \frac{\rho_n}{n} \bigl\{ l_c(\alpha)- l_c(
\alpha_0)- \bigl[ Q(\alpha)- Q(\alpha_0) \bigr] \bigr\}
\\
&&\qquad = \frac{1}{\tilde{\sigma}^2(\tilde{\alpha})}\frac{\rho_n}{n} \biggl\{ \bigl[(WY)^{\mathrm{T}}PA(
\tilde{\alpha})Y - L_n(\tilde{\alpha}) \bigr] -\frac{\tilde{\sigma
}^2(\tilde{\alpha})-\bar{\sigma}^2(\tilde{\alpha})} {
\bar{\sigma}^2(\tilde{\alpha})}L_n(
\tilde{\alpha}) \biggr\}
\\
&&\quad\qquad{}\times (\alpha-\alpha_0),
\end{eqnarray*}
where $\tilde{\alpha}$ lies between $\alpha$ and $\alpha_0$, and
$
L_n(\tilde{\alpha})= E[(WY)^{\mathrm{T}} P A(\tilde{\alpha})Y]$.
By the same arguments as in the proof of Theorem~\ref{teo1}, we have
$\tilde{\sigma}^2(\tilde{\alpha})
-\bar{\sigma}^2(\tilde{\alpha})=o_P(1) $ for any $\tilde{\alpha}$
on $\Delta$,
and $\bar{\sigma}^2(\alpha)$ is uniformly bounded away from zero on
$\Delta$.
So, $\tilde{\sigma}^2(\alpha)$ is uniformly bounded away from zero in
probability. This together with Lemmas~\ref{lem4}(5),~\ref{lem4}(6),~\ref{lem11}
and Chebyshev's
inequality lead to
\[
\frac{\rho_n}{n} \bigl\{ l_c(\alpha)- l_c(
\alpha_0)- \bigl[ Q(\alpha)- Q(\alpha_0) \bigr] \bigr
\}=o_P(1)\qquad\mbox{uniformly on } \Delta.
\]

The uniqueness condition of $\alpha_0$ can be obtained by Lemma~\ref{lem4},
Lemma~\ref{lem11},
and the same arguments as in the proof of Theorem~\ref{teo1}.
\end{pf*}

\begin{pf*}{Proof of Theorem~\ref{teo5}}
By Taylor's expansion, we have that
\[
0=\frac{\partial l_c(\hat{\alpha})}{\partial\alpha} =\frac{\partial
l_c(\alpha_0)}{\partial\alpha} +\frac{\partial^2 l_c(\tilde{\alpha
})}{\partial\alpha^2} (\hat{\alpha} -
\alpha_0),
\]
where $\tilde{\alpha}$ lies between $\hat{\alpha}$ and $\alpha_0$,
and thus
converges to $\alpha_0$ in probability by Theorem~\ref{teo4}. So, the asymptotic
distribution of $\hat{\alpha}$ can be obtained by proving that
\[
-\frac{\rho_n}{n} \frac{\partial^2 l_c(\tilde{\alpha})}{\partial\alpha
^2} \stackrel{P} {\longrightarrow}
\sigma_1^2 \quad\mbox{and}\quad\sqrt{ \frac{\rho_n}{n}}
\frac{\partial l_c(\alpha_0)}{\partial\alpha} \stackrel{D} {\longrightarrow}N \bigl(0, \sigma_2^2/
\sigma_0^4 \bigr),
\]
when\vspace*{-2pt} $\rho_n\rightarrow\infty$, where
$\sigma_1^2=\frac{1}{\sigma_0^2}\lim_{n\rightarrow\infty}
\frac{\rho_n}{n} E[(G\mathbf{m} -SG\mathbf{m})^{\mathrm
{T}}(G\mathbf
{m}-SG\mathbf{m})]$ and
$\sigma_2^2= \sigma_0^4\sigma_1^2$.

As we have, by $A(\alpha)A^{-1}=I_n + (\alpha_0 - \alpha)G$, Lemma
\ref{lem11} and
Chebyshev's inequality, that
$\frac{\rho_n}{n}(WY)^{\mathrm{T}}PWY =O_P(1)$
and
$\frac{\rho_n}{n}(WY)^{\mathrm{T}}P A(\alpha)Y = O_P(1)$,
so, when $\rho_n \rightarrow\infty$,
\[
\frac{\rho_n}{n}\frac{\partial^2 l_c(\alpha)}{\partial\alpha^2} = -\frac
{1}{\tilde{\sigma}^2(\alpha) }\cdot
\frac{\rho_n}{n}(WY)^{\mathrm{T}} P WY -\frac{\rho_n}{n} \operatorname{tr}
\bigl( \bigl[WA^{-1}(\alpha) \bigr]^2 \bigr)
+o_P(1).
\]
This together with Lemmas~\ref{lem6}(1),~\ref{lem8}(1) lead to
$
\tilde{\sigma}^2(\alpha)
=\sigma_0^2 + o_P(1)
$
for any $\alpha\in\Delta$ when $\rho_n\rightarrow\infty$. Therefore,
by the mean value theorem, conditions (5)--(6) and
$\tilde{\alpha}\stackrel{P}{\longrightarrow}\alpha_0$, we have
$
\frac{\rho_n}{n} \{
\frac{\partial^2 l_c(\tilde{\alpha})}{\partial\alpha^2}
-
\frac{\partial^2 l_c(\alpha_0)}{\partial\alpha^2}
\}=o_P(1)$.

It\vspace*{-1pt} follows, from $\tilde{\sigma}^2 (\alpha_0)\stackrel
{P}{\longrightarrow}\sigma_0^2$,
Lemma~\ref{lem11},
Chebyshev's inequality and the row sums of $G$ being uniform order
$O(1/\sqrt{\rho_n})$, that
$-\frac{\rho_n}{n}
\frac{\partial^2 l_c(\alpha_0)}{\partial\alpha^2}
\stackrel{P}{\longrightarrow}\sigma_1^2$.

In\vspace*{-1pt} the following, we will establish the asymptotic distribution of
$\sqrt{\frac{\rho_n}{n}} \frac{\partial l_c(\alpha_0)}{\partial
\alpha}$.

By\vspace*{-1pt} Lemmas~\ref{lem10}(2) and~\ref{lem11}(3), it is easy to see
$\sqrt{\frac{\rho_n}{n}}(G\mathbf{m})^{\mathrm{T}} P \mathbf{m}
=o_P(1)$ and
$\sqrt{\frac{\rho_n}{n}}(G\bolds{\varepsilon})^{\mathrm{T}}P\mathbf
{m} =o_P(1)$ when
$nh^6\rightarrow0$ and
$h^2 \log n\rightarrow0$. By straightforward calculations and Lemmas~\ref{lem6}(1),
\ref{lem8}(1),~\ref{lem11}(5) and~\ref{lem11}(7), we have the first-order derivative of
$\sqrt{\frac{\rho_n}{n}} l_c(\alpha)$ at $\alpha_0$ is
\[
\frac{1}{\tilde{\sigma}^2(\alpha_0)} \sqrt{\frac{\rho_n}{n}} \biggl\{ (G\mathbf{m} - SG
\mathbf{m})^{\mathrm{T}}\bolds{\varepsilon} + \bolds{\varepsilon
}^{\mathrm{T}} \biggl[G -\frac{1}{n}\operatorname{tr}(G)I_n
\biggr]\bolds{ \varepsilon} \biggr\} +o_P(1).
\]
By Lemma~\ref{lem12}, we have
\[
\sigma_{qn}^{-1} \biggl\{(G\mathbf{m} - SG
\mathbf{m})^{\mathrm
{T}}\bolds{\varepsilon} + \bolds{\varepsilon}^{\mathrm{T}}
\biggl[G^{\mathrm{T}}-\frac
{1}{n}\operatorname{tr}(G)I_n
\biggr]\bolds{\varepsilon} \biggr\} \stackrel{D} {\longrightarrow}N(0,1),
\]
where
$\sigma_{qn}^2 = \operatorname{var} \{(G\mathbf{m} - SG\mathbf
{m})^{\mathrm{T}}\bolds{\varepsilon}
+ \bolds{\varepsilon}^{\mathrm{T}}[G -\frac{1}{n}\operatorname
{tr}(G)I_n]\bolds{\varepsilon} \}$.
So, by $\frac{\rho_n}{n}\sigma_{qn}^2\rightarrow\sigma_2^2$
and $\tilde{\sigma}^2(\alpha_0)\stackrel{P}{\longrightarrow}\sigma
_0^2$, we have
$
\sqrt{\frac{n}{\rho_n}}(\hat{\alpha}-\alpha_0)
\stackrel{D}{\longrightarrow}N (0, \sigma_0^2\lambda_4^{-1} )$.
\end{pf*}

\begin{pf*}{Proof of Theorem~\ref{teo6}}
By straightforward calculations, Lemmas~\ref{lem6}(1),
\ref{lem8}(1),~\ref{lem11}, Chebyshev's inequality
and Theorem~\ref{teo5}, we have
$
\sqrt{n}(\hat{\sigma}^2-\sigma_0^2)
=\frac{1}{\sqrt{n}}\sum_{i=1}^n(\varepsilon_i^2-\sigma_0^2)+o_P(1)
$
when $\rho_n\rightarrow\infty$. This together with the central limit theorem
lead to Theorem~\ref{teo6}.
\end{pf*}

\begin{pf*}{Proof of Theorem~\ref{teo7}}
Theorem~\ref{teo7} can be obtained by using the same
arguments as in the proof of Theorem~\ref{teo3}, except that here
\begin{eqnarray*}
J_{n2}&=&\sqrt{f(s)}\sqrt{\frac{n h_1^2}{\rho_n}}(
\alpha_0-\hat{\alpha}) (I_p, \mathbf{0}_{p\times2p})
\bigl(n^{-1}H_1^{-1}\mathcal{X}_1^{\mathrm
{T}}
\mathcal{W}_1\mathcal{X}_1 H_1^{-1}
\bigr)^{-1}
\\
&&{}\times \frac{\sqrt{\rho_n }}{n}H_1^{-1}
\mathcal{X}_1^{\mathrm{T}}\mathcal{W}_1 G(\mathbf{m} +
\bolds{\varepsilon}).
\end{eqnarray*}
By Lemma~\ref{lem2}(1), Markov's inequality, the row sums of the matrix $G$ having
uniform order $O(1/\sqrt{\rho_n})$ and condition~(4), we have
\[
(I_p, \mathbf{0}_{p\times
2p}) \bigl(n^{-1}H_1^{-1}
\mathcal{X}_1^{\mathrm{T}}\mathcal{W}_1\mathcal
{X}_1 H_1^{-1} \bigr)^{-1}
\frac{\sqrt{\rho_n}}{ n} H_1^{-1}\mathcal{X}_1^{\mathrm
{T}}
\mathcal{W}_1 G(\mathbf{m}+\bolds{\varepsilon})=O_P(1).
\]
Furthermore, it can be seen from the proof of Theorem~\ref{teo5} and Lemma~\ref{lem10}(1) that
when $nh_1^6=O(1)$ and $h/h_1\rightarrow0$,
$\sqrt{\frac{n h_1^2}{\rho_n}}(\hat{\alpha}-\alpha)\stackrel
{P}{\longrightarrow}0$. So,
$J_{n2}=o_P(1)$.
\end{pf*}
\end{appendix}


%
\begin{supplement}
\stitle{Detailed proofs of lemmas and theorems}
\slink[doi]{10.1214/13-AOS1201SUPP} 
\sdatatype{.pdf}
\sfilename{aos1201\_supp.pdf}
\sdescription{We provide the detailed proofs of the lemmas and theorems.}
\end{supplement}


%

\printaddresses

\end{document}